\documentclass[a4paper,11pt]{amsart}
\usepackage{amsthm}
\usepackage{amssymb}
\usepackage[totalwidth=480pt,totalheight=680pt]{geometry}
\def\Arg{\operatorname{Arg}}

   \def\pdiff#1#2{\dfrac{\partial#1}{\partial#2}}

\def\d{\operatorname{d}\!}

\def\thei{{\bold i}}

\def\modulus#1{\left|#1\right|}
\def\norm#1{\left\|#1\right\|}
      \def\normp#1_#2{\norm{#1}_{#2}}

\def\Set#1{\left\{\,#1\,\right\}}

\let\Q\Rational

\let\R\Real

\let\sphere\Sphere

\def\map#1#2#3{#1~\colon~#2~\to~#3}
\input epsf.tex
\def\cal{\mathcal}

\def\DEF{\stackrel{\scriptstyle\text{def}}{:=\!=}}

   \def\patrix#1{\begin{pmatrix}#1\end{pmatrix}}

\title[Off-center Reflections]
{Off-center Reflections: Caustics and Chaos}
\author[T.~K.~K.~Au]{Thomas Kwok-keung Au}
\address{Dept of Math, The Chinese University of Hong Kong}
\email{thomasau@cuhk.edu.hk}
\thanks{The first named author is partially supported by 
UGC Earmarked Grant and CUHK Direct Grant. The second named 
author would like to thank the faculty and staff of
the math department of CUHK for their hospitality during his visits.}
\author[X.~S.~Lin]{Xiao-song Lin}
\address{Dept of Math, University of California, Riverside}
\email{xl@math.ucr.edu}

\newtheorem{theorem}{Theorem}[section]
\newtheorem{proposition}[theorem]{Proposition}
\newtheorem{lemma}[theorem]{Lemma}

\newtheorem{conjecture}[theorem]{Conjecture}

\begin{document}
\pagestyle{empty}
\pagestyle{headings}
\setcounter{page}{1}
\maketitle

\setlength\topmargin{0mm}
\setlength\textheight{230mm}
\setlength\oddsidemargin{0mm}
\setlength\evensidemargin\oddsidemargin
\setlength\textwidth{160mm}

\thispagestyle{empty}
\setlength\baselineskip{24pt}
\setlength\parindent{3.25em}
\setlength\parskip{1ex}

\section{Introduction}
\label{Sec-Intro}
We study the properties of a particular one-parameter family of
circle maps 
called {\sl off-center reflections} defined in \S\ref{Sec-basic}.
This map, in its 2-dimensional version, is
first introduced in an open problem by S.~T.~Yau,
\cite[problem~21]{Y},
who suggests a cross study of the dynamics and geometry.
In this article, we attempt to explore the possible link between the
dynamics
of this family of circle maps and their {\sl caustics}.
Although our study has not much
contents in differential geometry as Yau expected,
it reveals some interesting phenomena.
For example, we observe and partially 
prove that in a certain generic range of the parameter,
the caustics have {\it exactly} 4 cusp points for odd iterations;
whereas for even iterations, each caustic is a curve
tangential to the circle at {\it exactly} four points.  This may not
be the best result that one could state about the dynamics and the
geometry of the map; nevertheless, we still put it forward in the
hope that our study may invite better understanding to the subject.
The off-center reflection also bares several interesting analytic
forms.  It is a Blaschke product restricted to the circle.
Moreover, it has an infinite series expression highlighting that it
is a perturbation of rotation on the circle.
Since the work of Arnold, \cite{A2}, a standard type of perturbations
has attracted much interests in mathematics and
physics communities, \cite{BBJ,Di,Z}.  This standard type is exactly
a reduction of the series of the off-center reflection.
This adds more flavor to our study.
 
This family of off-center reflections plays an interesting
role in the space of circle maps.
With the parameter $r$ going from 0 to 1, it carries the
initial antipodal map to the terminal doubling map, which
provides a particularly nice way of deforming a simple
dynamics to a chaotic one (see the
asymptotic orbit diagram in \S\ref{Sec-details}). We
expect that some notion of stability about the
aforementioned cusp points and tangent points on the caustics
might emerge from further study of this family of circle maps.
We also hope to explore more systematically the dynamics of
general circle maps through the method of symplectic and
contact geometry in a forthcoming work.
  
We will begin in \S\ref{Sec-basic} with the definition and
some analytic properties of
the map.  Then the study is divided into three parts.
In the first part (\S\ref{Sec-caustics0}), our attention is
given to the caustics of the map and its iterations.
Our results concerning the caustics of odd iterations are more
conclusive.  Following a symplectic and contact geometry
interpretation developed by Arnold \cite{A1},
we discover the \lq\lq generating function'' for the corresponding
Lagrangian embedding.  This in turns provides explicit formulae for
the {\sl orthotomics} of the caustics.
Moreover,
these orthotomics are always smooth simple closed curves
(and hence the caustic of an off-center reflection has at
least 4 cusp points following the classical 4 vertex
theorem), and they are convex
in a certain range of the parameter.  The method fails for even
iterations.  Nevertheless, explicit computations
still provide reasonable support for our prediction.

In the second part (\S\ref{Sec-experiments}), the main
observations are presented together with graphical illustrations.
There are also theoretical support for them.  For
example, we have partial result that the caustics is stable with
$r\leq 1/3$.  The tedious
computations are separated into
\S\ref{Sec-details} for detailed reading.
 
In the third part (\S\ref{Sec-modelock}), 
we study the phenomenon of mode-locking behavior
for this particular family of circle maps and the width of resonance
zone is estimated.  This is an attempt to
understand the iterations of the map.  This family extends a
class
of examples, which exhibits the same behavior, studied by Arnold
and others.  The mode-locking phenomenon of the off-center
reflections and its ``complex conjugates'' are totally different.
Moreover, $r=1/3$ is the first value that this behavior undergoes a
structural change.  This probably is not simply a coincidence with the
bifurcation values of cusps.  There is room for investigation in this
aspect.

\section{Off-Center Reflection}
\label{Sec-basic}
An {\sl off-center reflection} is a map $\sphere^1\rightarrow \sphere^1$ defined
as follows: Pick a point, say $(r,0)$ in the interior of the unit disk $D^2$. For any $z\in
\partial  D^2=\sphere^1$, emit a ray from $(r,0)$ to $z$. The ray will
be reflected at $z$
with $\sphere^1$ as the curve of reflection and the reflected ray will hit
another point $R_r(z)$ on $\sphere^1$. The map $z\mapsto R_r(z):\sphere^1\rightarrow \sphere^1$ is
what we call an off-center reflection.   The action of the map is
shown in the figure where a point on $\sphere^1$ is represented by
$\phi\mod 2\pi$ and $\alpha$ denotes the incident angle.
\begin{center}
\mbox{\epsfysize=50mm \epsfbox{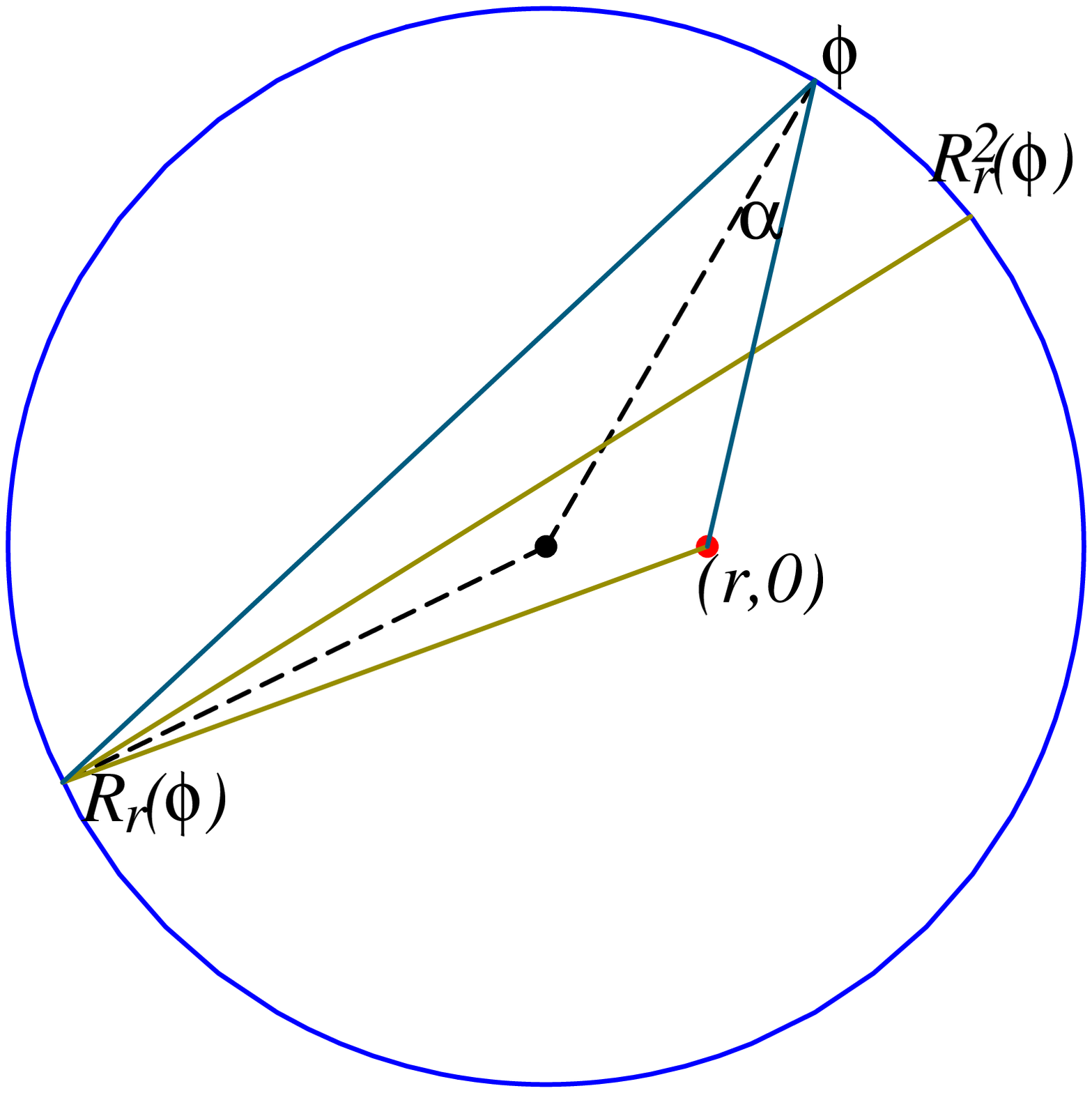}}
\end{center}

The iteration of this map is a little uncommon at first
sight because it is different from the usual successive
reflection in a curved mirror.  However,
maps similar to it
have been a center of discussion in circle dynamics.
Let us first establish two analytic
expressions for the map $R_r:\sphere^1\rightarrow \sphere^1$ where
$0\leq r < 1$.  By them, we may see $R_r$ as a real function
in terms of an infinite series as well as a complex function
restricted to the circle.
They are
\begin{gather*}
R_r(\phi)=\phi+\pi-2\alpha\qquad \text{mod}\,\,2\pi, \\
\alpha =\alpha(\phi) \DEF \Arg\left(\cos{\phi}-r+\thei\sin{\phi}\right)-\phi,
\end{gather*}
where $\Arg$ is the principal argument taking values in
$(-\pi,\pi]$.
\begin{lemma}
The angle of incident $\alpha$ has the following Fourier sine series,
$$
\alpha(\phi) = \Arg\left(\cos{\phi}-r+\thei\sin{\phi}\right)-\phi
=\sum_{k=1}^\infty\frac{r^k}{k}\sin(k\phi).
$$
Thus,
\begin{equation}\label{EqnRrSeries}
R_r(\phi) =
\phi+\pi-2\sum_{k=1}^\infty
\frac{r^k}{k}\sin(k\phi)\qquad\text{\rm mod}\,\,2\pi.
\end{equation}
\end{lemma}
\begin{proof}
It is clear that both $\alpha$ and $\partial\alpha/\partial r$ are odd
functions in $\phi$.  So, they both have Fourier sine series expansions.  Let
$\alpha = \sum a_k(r)\sin(k\phi)$.  Then, the
coefficients of the series for $\partial\alpha/\partial r$ are given by
\begin{align*}
\frac{\partial a_k}{\partial r} &= \frac{2}{\pi}\int_0^{\pi} \frac{ \sin
\phi\sin(k\phi) }{1 - 2r\cos
\phi + r^2}
\d \phi
\\ &= \frac{1}{\pi} \int_0^{\pi} \left(\frac{ \cos(k-1)\phi }{1 - 2r\cos\phi + r^2} -
\frac{\cos(k+1)\phi }{1 - 2r\cos\phi + r^2}\right) \d \phi \\
&= \frac{1}{\pi}\left( \frac{\pi r^{k-1}}{1-r^2} - 
\frac{\pi r^{k+1}}{1-r^2}\right) = r^{k-1}.
\end{align*}
Hence, those coefficients $a_k$ for the series of $\alpha$ must be
$r^k/k$.
\end{proof}

In this paper, we often omit the $\text{mod } 2\pi$ when it is clear in 
context.  Formula~(\ref{EqnRrSeries}) of $R_r$ without the
modulo~$2\pi$ is exactly the lifting of $R_r$ to a function from
$\R$ to $\R$ taking $0$ to $\pi$.

By playing with the argument of a complex number, we 
get another expression for the map $R_r$.  This shows that the
off-center reflection is a special case of the so-called Blaschke
product.  This map is not extendible to the hyperbolic disk.
\begin{lemma}
For $|z|\leq 1$, we have a complex function
\begin{equation}\label{EqnRrCpx}
z \mapsto -z^2\,\frac{1-rz}{z-r},
\end{equation}
whose restriction to the boundary of the disk, $\modulus{z}=1$,
is the off-center reflection $R_r$.
\end{lemma}

Therefore, this function
$R_r(\phi)$ is harmonic when $(r,\phi)$ are treated as polar 
coordinates as it is the argument of the analytic function
$-z(1-z)^2$.  We do not know whether this coupling between 
the parameters has additional physical or geometrical
implications.

By changing a sign of the off-center reflection, we have another map
$\overline{R}_r$ defined by
$$
\overline{R}_r ~:~\phi \mapsto \phi+\pi+2\sum_{k=1}^\infty
\frac{r^k}{k}\sin(k\phi)\qquad\text{mod}\,\,2\pi.
$$
Geometrically, it is the mirror image of $R_r(\phi)$ reflected 
by the diameter joining $\phi$ to $\phi+\pi$.   
This map can be extended to the unit disk, namely,
$$
z \mapsto -\frac{z-r}{1-rz}.
$$
Therefore, it defines a map in $\text{PSL}(2,\R)$, the isometry group of the
hyperbolic disk.  The dynamics of $R_r$ and this ``conjugated'' map
$\overline{R}_r$ are completely different 
(see \cite{H2} and \S\ref{Sec-width}).

For sufficiently small $r$, $R_r$ behaves very similarly to $R_0$,
i.e. the antipodal map.  In fact, when $r<1/3$, it is in the same
component of
$R_0$ in the group of orientation preserving diffeomorphisms of
$\sphere^1$.  However, $R_{1/3}$ is only a
homeomorphism on $\sphere^1$ and $R_r$ a degree~1 map when $1/3 < r
< 1$.  These can be easily concluded from the derivatives of $R_r$,
which will also be useful later,
\begin{align*}
R'_r(\phi) &= \dfrac{1-4r\cos\phi + 3r^2}{1-2r\cos\phi+r^2}, \\
R''_r(\phi) &= \frac{2r(1-r^2)\sin\phi}{(1-2r\cos\phi+r^2)^2}, \\
R_r^{(3)}(\phi) &= \frac{2r(1-r^2) \left[(1+r^2)\cos\phi -
2r(1+\sin^2\phi)\right]}{(1-2r\cos\phi+r^2)^3}.
\end{align*}

In \S\ref{Sec-details}, we will give more information about the
fixed point and other special points of $R_r$.  More
dynamical properties such as periodic cycles and whether
they are attracting are discussed in \cite{Au}.

\section{Caustics}
\label{Sec-caustics0}

\subsection{Caustic of Off-center Reflection}
\label{Sec-caustics}
In this section, we will discuss the cusp phenomenon of the
caustic of $R_r$.  Classical examples of caustics are the locus of
focal points with respect to a point on a surface and the focal curve of a
convex plane curve.  Corresponding to these caustics, there
are the famous Geometric Theorem (Conjecture) of Jacobi and Four-vertex
Theorem.
There are many interesting at-least-four results, see
\cite{A1,A5,T1,T2}.  The caustic of off-center reflection
provides another one.  The conjugate locus of a point on a flat
flying disc is, at degenerate situation, the caustic of the
off-center reflection.  It should be remarked that in the series of
papers \cite{BGG1,BG1,BG2,GK}, the authors analyzed the
singularities of the caustics produced by a point light source when
it is reflected in a codimension~1 ``mirror'' in $\R^2$ and $\R^3$.
Their emphasis though is on the ``source genericity'': whether the
caustics could be made generic by moving the source.  See also
\cite{BGG2}.

For a circle map $\map{f}{\sphere^1}{\sphere^1}$,
the family of lines joining $\phi$ to $f(\phi)$ is
\begin{align*}
F(\phi,x,y)&=(\sin f(\phi)-\sin\phi)(x-\cos\phi)-(\cos
f(\phi)-\cos\phi)(y-\sin\phi) \\
&=(\sin f(\phi)-\sin\phi)x-(\cos
f(\phi)-\cos\phi)y-\sin\left(f(\phi)-\phi\right).
\end{align*}
The {\sl caustic} of the map $f$ is defined to be the envelope of these
lines.  Thus, it is given by the equations 
$\pdiff{F}{\phi}(\phi,x,y) = 0 = F(\phi,x,y)$, that is,
$$
\begin{pmatrix}
\sin f(\phi) - \sin\phi & -\cos f(\phi) + \cos\phi \\
f'(\phi)\cos f(\phi) - \cos\phi & f'(\phi)\sin f(\phi) -
\sin\phi
\end{pmatrix}
\begin{pmatrix} x \\ y \end{pmatrix} =
\begin{pmatrix}
\sin(f(\phi)-\phi) \\
\left(f'(\phi)-1\right)\cos(f(\phi)-\phi)
\end{pmatrix}
$$
Solving for $x, y$, we obtain a parameterization of the
caustic
\begin{equation}
\left\{\begin{aligned}
x(\phi) &=\frac{f'(\phi)\cos\phi+\cos f(\phi)}{1+f'(\phi)}\\
y(\phi) &=\frac{f'(\phi)\sin\phi+\sin f(\phi)}{1+f'(\phi)}.
\end{aligned}\right.\label{Eqncausticf}
\end{equation}
The tangent direction, which is degenerated at cusp points,
of the caustic is given by
\begin{equation}
\left\{\begin{aligned}
x'(\phi)&=\frac{f''(\phi)(\cos\phi-\cos f(\phi))-
f'(\phi)(1+f'(\phi))(\sin\phi+\sin
f(\phi))}{(1+f'(\phi))^2}\\
y'(\phi)&=\frac{f''(\phi)(\sin\phi-\sin f(\phi))+
f'(\phi)(1+f'(\phi))(\cos\phi+\cos f(\phi))}
{(1+f'(\phi))^2}.
\end{aligned}\right.\label{Eqncausticftgt}
\end{equation}

The caustic (\ref{Eqncausticf}) of the off-center reflection
may run to infinity since
$
1+R'_r(\phi) 
$
may be equal to zero.
In fact, it is so if and only if
$r\geq 1/2$.  It should be more appropriate to define the caustic as
the envelop of the geodesic normal field on the sphere.  After the
stereographic projection, it does not matter whether the caustic is
defined on the plane or the sphere as the local
properties remain unchanged 
(Darboux Theorem of symplectic structure).  As we will see, the
local property of the caustic of $R_r$ can be understood by direct
computation.
\begin{theorem}\label{ThmExact4Rr}
For all $0<r<1$, there are exactly 4 cusp points 
on the caustic of $R_r$. Two of them correspond to the $R_r$-orbit,
$\{0,\pi\}$, of period 2.
\end{theorem}
\begin{proof}
The derivatives of $x$ and $y$ can be expressed
completely in terms of
$r$ and
$\phi$, namely,
\begin{align*}
x'(\phi) &= \frac{6r^2(-\cos\phi+r\cos(2\phi))(r-\cos\phi)\sin\phi}{\left( -1 - 2\,{r^2} +
3\,r\,\cos\phi \right)^2} \\
y'(\phi) &= \frac{6r^2(-1+2r\cos\phi)(r-\cos\phi)\sin^2\phi}{\left( -1 - 2\,{r^2} +
3\,r\,\cos\phi \right)^2}.
\end{align*}
The common solutions for $x'(\phi) = 0 = y'(\phi)$ are 
$\phi = 0, \pi$ and two values of $\phi$ with $\cos\phi = r$.
Clearly, $0$ and $\pi$ are zeros of $x'$ of first order and
of $y'$ of second order, thus,
these are semicubical cusps.  If $\cos\phi = r$, after
further differentiation and evaluation at the point, one has
$$
\begin{aligned}
x''(\phi) &= -12r^3; \\
y''(\phi) &= \frac{6r^2(2r^2-1)}{\sqrt{1-r^2}}; \\
\end{aligned}
\qquad\qquad
\begin{aligned}
x'''(\phi) &= \frac{12(5r^4+r^2)}{\sqrt{1-r^2}}; \\
y'''(\phi) &= \frac{-6r^3(10r^2-3)}{1-r^2}.
\end{aligned}
$$
Thus, 
$
x''y''' - x'''y'' = \frac{72r^4}{1-r^2}\ne 0.
$
Therefore, there are also semicubical cusps at those values of
$\phi$ with $\cos\phi = r$.
\end{proof}
Here are two pictures of the caustics of $R_r$, one for $r < 1/2$ and
the other $r>1/2$.  Since the second
one runs to infinity, it is drawn with ``compressed'' scale
where a circle of radius~$>1$ represents the point of infinity and
the caustics has a self-intersection there.
\begin{center}
\parbox{75mm}{\begin{center}
\mbox{\epsfysize=43mm \epsfbox{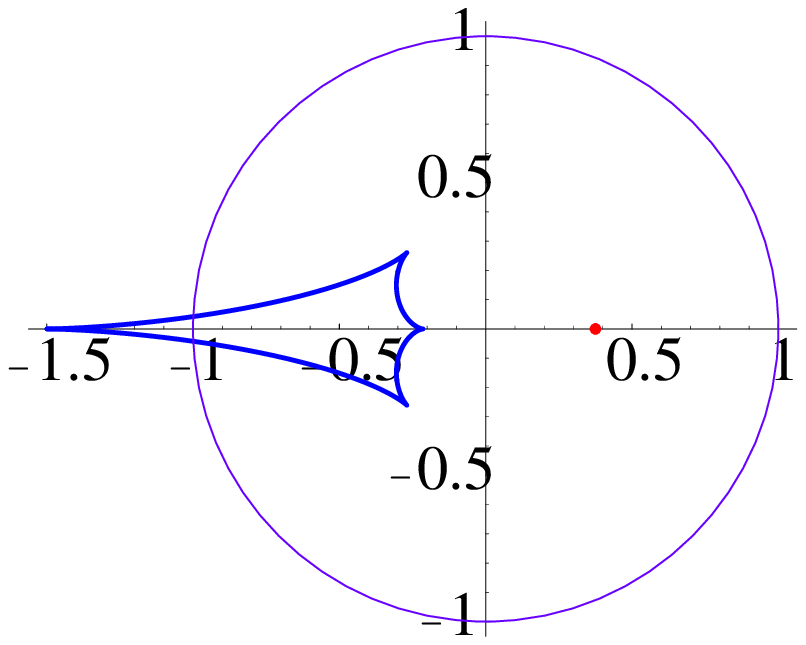}} \\
{\footnotesize caustic of $R_r$ with $r=0.375$}
\end{center}} \hfil
\parbox{75mm}{\begin{center}
\mbox{\epsfysize=43mm \epsfbox{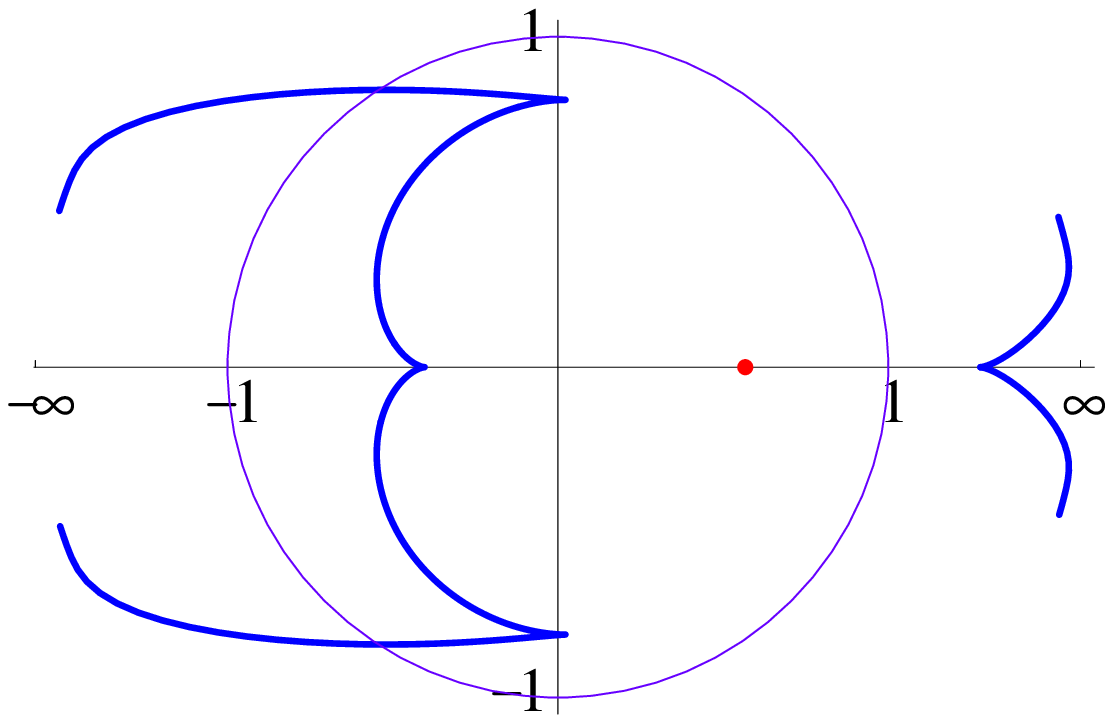}} \\
{\footnotesize caustic of $R_r$ with $r\approx 1/\sqrt{2}$}
\end{center}}
\end{center}
%

\subsection{Symplectic and contact geometry interpretation}
\label{Sec-symplectic}
The explicit computation in the previous section has its
advantages and shortcomings.  On the one hand, it gives a very
exact count of the number of cusp points; on the other
hand, it is too specific and also very complicated to apply 
especially when iterations are considered.  In this section, we
present the symplectic method which may work in general,
though the result obtained is not as specific as before.  The
terminology of \cite{A1} is followed.

Denote the coordinates of the unit cotangent bundle
$ST^*(\R^2)$ by $(p_x,p_y,x,y)$ where $(x,y)\in\R^2$ and
$p_x^2+p_y^2=1$. This bundle
is a contact 3-manifold with the contact
1-form $p_x\,dx+p_y\,dy$,
and the cotangent manifold $T^*(\R^2)$ is symplectic with the 
symplectic 2-form $d\,(p_x\,dx+p_y\,dy)$.

With our notation of $\phi$ and $\alpha$ in \S\ref{Sec-basic},
the vector from $\phi$
to $R_r(\phi)$ is never zero and the unit vector in that direction
$(p_x,p_y)$ is given by 
$$p_x=\cos(\phi+\pi-\alpha),\qquad p_y=\sin(\phi+\pi-\alpha).$$
Then
$$
\patrix{\phi \\ S} \mapsto
\patrix{p_x \\ p_y \\ x \\ y} =
\patrix{\cos(\phi+\pi-\alpha) \\ \sin(\phi+\pi-\alpha) \\
\cos\phi+S\cos(\phi+\pi-\alpha) \\
\sin\phi+S\sin(\phi+\pi-\alpha)}
$$
defines a map
$L:\sphere^1\times\R^1\longrightarrow T^*(\R^2).$
This map may be thought of as a flow (in parameter $S$) of
unit speed in the direction of the reflection lines,
starting with the round circle $S=0$.  This is a case of
what Arnold called \lq\lq Legendrian collapsing" \cite{A1}. 

Let $p:T^*(\R^2)\rightarrow\R^2$ be the canonical
projection.  The Jacobian of the map $p\circ L$ is
\begin{align*}
J(p\circ L) &=\text{det}\pmatrix \cos(\phi+\pi-\alpha) & 
	 -\sin\phi-S(1-\alpha')\sin(\phi+\pi-\alpha)\\
	\sin(\phi+\pi-\alpha) &
	  \cos\phi+S(1-\alpha')\cos (\phi+\pi-\alpha)
	\endpmatrix \\
	    &=\cos(\phi+\pi-\alpha)\cos\phi+\sin(\phi+\pi-\alpha)\sin\phi
		+S(1-\alpha')\\
	    &=-\cos\alpha+S(1-\alpha')
\end{align*}
So the equation for the {\sl critical curve} on 
$\sphere^1\times\R^1$ is
$$S=\frac{\cos\alpha}{1-\alpha'}.$$

\begin{proposition}
The critical curve, when mapped to the $(x,y)$-plane, agrees with the 
caustic of $R_r$.
\end{proposition}
\begin{proof}
Notice that since $R_r(\phi)=\phi+\pi-2\alpha$, we have
$$R_r'=1-2\alpha'\qquad\text{or}\qquad 1+R_r'=2(1-\alpha').$$
Therefore, 
\begin{align*}
x&=\cos\phi+\frac{\cos\alpha}{1-\alpha'}\cos(\phi+\pi-\alpha)\\
&=\frac{(1+R_r')\cos\phi+2\cos\alpha\cos(\phi+\pi-\alpha)}{1+R_r'}\\
&=\frac{R_r'\cos\phi+\cos(\phi+\pi-2\alpha)}{1+R_r'},
\end{align*}
which agrees with one of the equations (\ref{Eqncausticf}). 
We may get the expression for $y$ similarly.
\end{proof}

Furthermore, a straightforward calculation shows that
$$p_x\,dx+p_y\,dy=\sin\alpha\,d\phi+dS.$$
Thus, the image of $L$ is a Lagrange cylinder in $T^*(\R^2)$. 
Notice further that $\alpha=\alpha(\phi)$ is an odd 
function of $\phi$ and therefore
$$\int_{\sphere^1}\sin\alpha\,d\phi=0.$$
This implies that $p(L)$ is an exact Lagrange cylinder.
Take a function $S(\phi)$ along the circle by
$$S(\phi)=-\int_0^{\phi}\sin\alpha\,d\phi.$$
It defines a section in the Lagrange cylinder $p(L)$.
It is easy to see that $-S(\phi)$ is increasing 
for $0\leq\phi
\leq\pi$ and decreasing when $\pi\leq\phi\leq2\pi$. Therefore, 
the curve $C$ given by
\begin{align*}
&x=\cos\phi+S(\phi)\cos(\phi+\pi-\alpha)\\
&y=\sin\phi+S(\phi)\sin(\phi+\pi-\alpha)
\end{align*}
is quite likely to be a convex plane curve. 
Let us show that this is the case when $r\leq1/2$.

First we note that $C$  has a continuous normal field 
$(\cos(\phi+\pi-\alpha),\sin(\phi+\pi-\alpha))$. It is easy to compute
$${x'}^2+{y'}^2=4\sin^2\alpha+(\cos\alpha-S(1-\alpha'))^2.$$
Therefore, ${x'}^2+{y'}^2=0$ is possible only when $\phi=\pi$. 
But the number of zeros of ${x'}^2+{y'}^2$ should be even
(geometrically, because $C$ is co-oriented). 
Thus ${x'}^2+{y'}^2>0$ all the time
and $C$ is a smooth simple closed curve. 
We also point out that the curvature of $C$ is
$$\kappa=\frac{1-\alpha'}{\sqrt{{x'}^2+{y'}^2}}$$
which is non-negative when $r\leq 1/2$ 
and therefore $C$ is convex.
Now the family of reflection lines of $R_r$ is identical to
the family of normal lines of this convex curve $C$. 
Therefore, the caustic of 
$R_r$ has at least 4 cusp points \cite{A1}.

The function $S(\phi)$ should be thought of as the 
{\sl generating function} of the circle map $R_r(\phi)$. 
The curve $C$ is related to the {\sl orthotomic} of such
a reflection.  We will study such generating functions for
general circle maps in a forthcoming work.

\subsection{Iterations of Reflections}
\label{Sec-iterate}
For an integer $n$, we denote the iteration of $R_r$ by
$R^n_r=R_r\circ R^{n-1}_r:\sphere^1\rightarrow \sphere^1$.
The equations (\ref{Eqncausticf}) and (\ref{Eqncausticftgt}) 
in \S\ref{Sec-caustics}, with $f=R^n_r$, give a
parametrization of the caustic of $R_r^n$ and its tangent.

The cusps on caustics of $R_r^n$ are more intriging and
complicated than that of $R_r$.  There are fundamental differences between the
caustics when $n$ is odd or even.  To see this difference, we may consider the
trivial example that $r = 0$.  For any odd $n$, $R_0^n$ is the antipodal map and
its caustic is a point (a degenerated curve with cusp).  
However, for even $n$, we
have the identity map, so the caustic is defined by the family 
of tangents and it is
the circle itself  (a smooth curve).  It is expected that this 
cusp versus ``smooth'' situation remains for $r$ close to $0$.  

\begin{theorem}\label{ThmAtleast4Rodd}
For sufficiently small $r>0$, the caustic of $R_r^{2m+1}$ has at
least 4 cusp points.
\end{theorem}
\begin{proof}
To some extent, the symplectic method in \S\ref{Sec-symplectic} may 
be adopted for
$R_r^{2m+1}$.  It can be proved by induction on $m$ that
$$
R_r^{2m+1}(\phi) = \phi + \pi - 2\tilde{\alpha}_m(\phi)
$$
for some odd function $\tilde{\alpha}_m(\phi)$.
For example,
$$
\tilde\alpha_1(\phi)=\alpha(\phi)+\alpha(\phi+\pi-2\alpha(\phi))
+\alpha(\phi-2\alpha(\phi)-2\alpha(\phi+\pi-2\alpha(\phi))).
$$
Therefore, we still have an exact Lagrangian cylinder and a
sectional curve defined by
$$
\tilde S(\phi)=-\int_0^{\phi}\sin\tilde\alpha_m\,d\phi.
$$
For small $r$, it is a convex curve and hence there
are at least 4 cusp points on the caustic.
\end{proof}
\noindent{\em Remark.\/}
The argument fails for even iterations of
$R_r$, because the analog of $\tilde\alpha$ is not an odd function.
Furthermore, it is not clear about how small the range of
$r$ should be.  Yet, from experimental observation, there are
at least four cusps for any $r>0$ and there are exactly four for 
$0<r<1/3$.  For more, please see the discussion after
proposition~\ref{PropCusp0piRodd}.

\section{Experiments and Observations}
\label{Sec-experiments}
To get more accurate information about the caustics of
iterations, we have to rely on lengthy calculations.  Our
investigation is indeed partly theoretical and partly
experimental.  We will first describe some interesting
properties with illustrations.  The technical details of
justification are left to the interested reader in \S\ref{Sec-details}.
\subsection{Observations}
\label{Sec-observations}
We first put forward a conjecture
about the exact picture of the caustic when $R_r^n$ is still a
diffeomorphism.
Then we look at the bifurcation process of the structure of the
caustics when $r$ varies.
Finally, we compare the caustics for different $n$.
We will soon see that the 2-cycles of $R_r$ play a special role
(propositions~\ref{PropTgtReven} and~\ref{PropQuadRinfty}).  The
2-cycles are $\Set{0,\pi}$ and $\Set{\pm\phi_c}$ where $\phi_c \in
(0,\pi)$ and $\cos\phi_c = \dfrac{1-\sqrt{1+8r^2}}{4r}$.  We will
often refer to this notation.

\begin{conjecture}\label{ConjAllRr}
For $0 < r \leq 1/3$, the caustic of $R_r^{2m+1}$ is a C$^\infty$
curve with exactly four cusp singularities, with two of them 
occurring at $\phi=0, \pi$.  
On the other hand, the caustic of $R_r^{2m}$ is a differentiable curve;
C$^\infty$ everywhere except at exactly the four 2-periodic 
points of $R_r$, where the caustic is tangent to the unit circle.
\end{conjecture}
The conjecture is demonstrated by the following pictures,
which are produced by programming in Mathematica.  
The purple (thin) curve is the unit circle and the
blue (thick) curve is the caustic with the light source
at the red dot, $(r,0)$.
Besides these pictures, we also have
partial results that support our conjecture.

\medskip
{\parindent=0mm
\parbox{52mm}{\centering%
\mbox{\epsfxsize=50mm \epsfbox{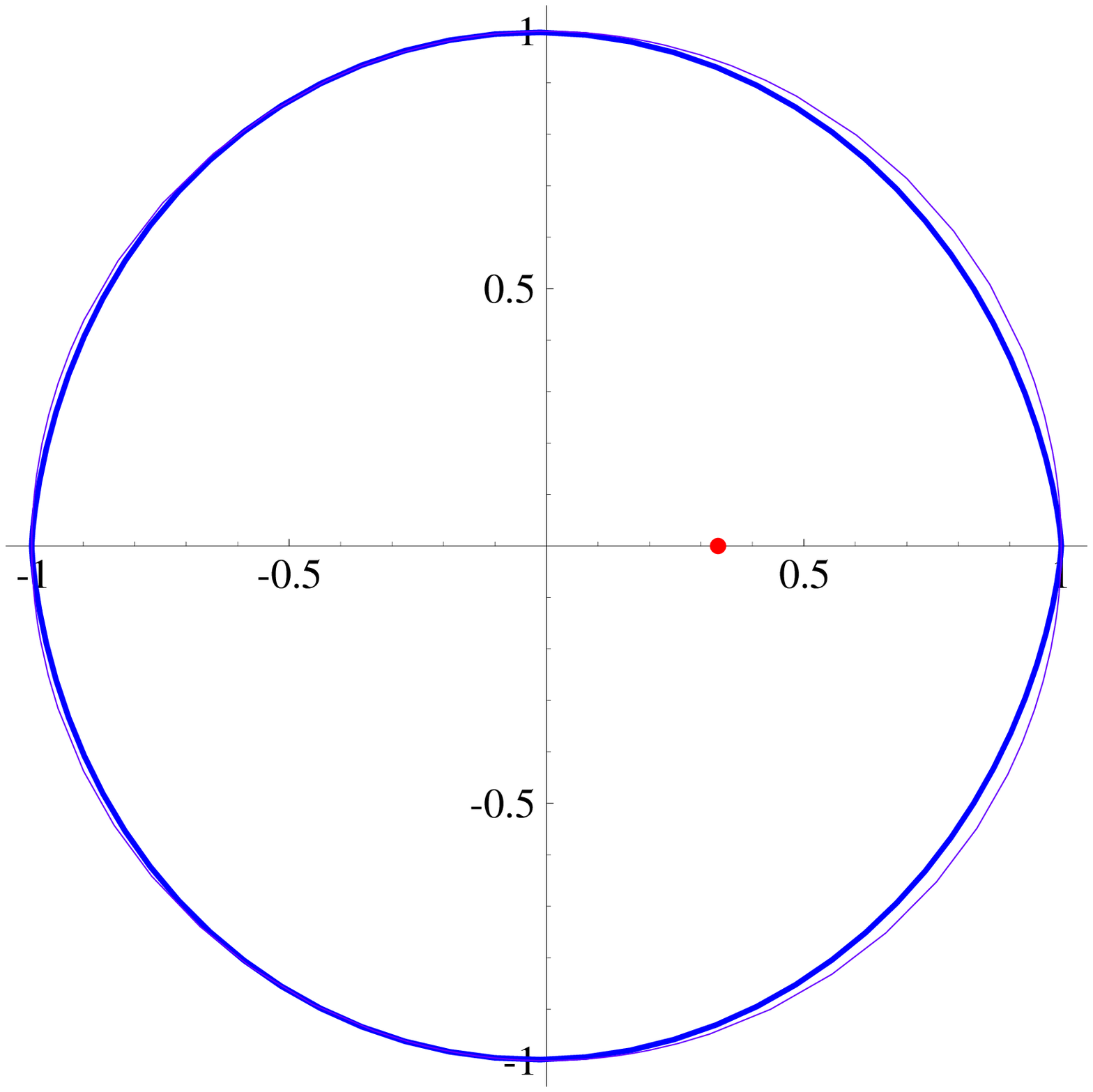}}\\
{\tiny caustic for $n=2$, $r=1/3$.}
} \hfill
\parbox{52mm}{\centering%
\mbox{\epsfxsize=50mm \epsfbox{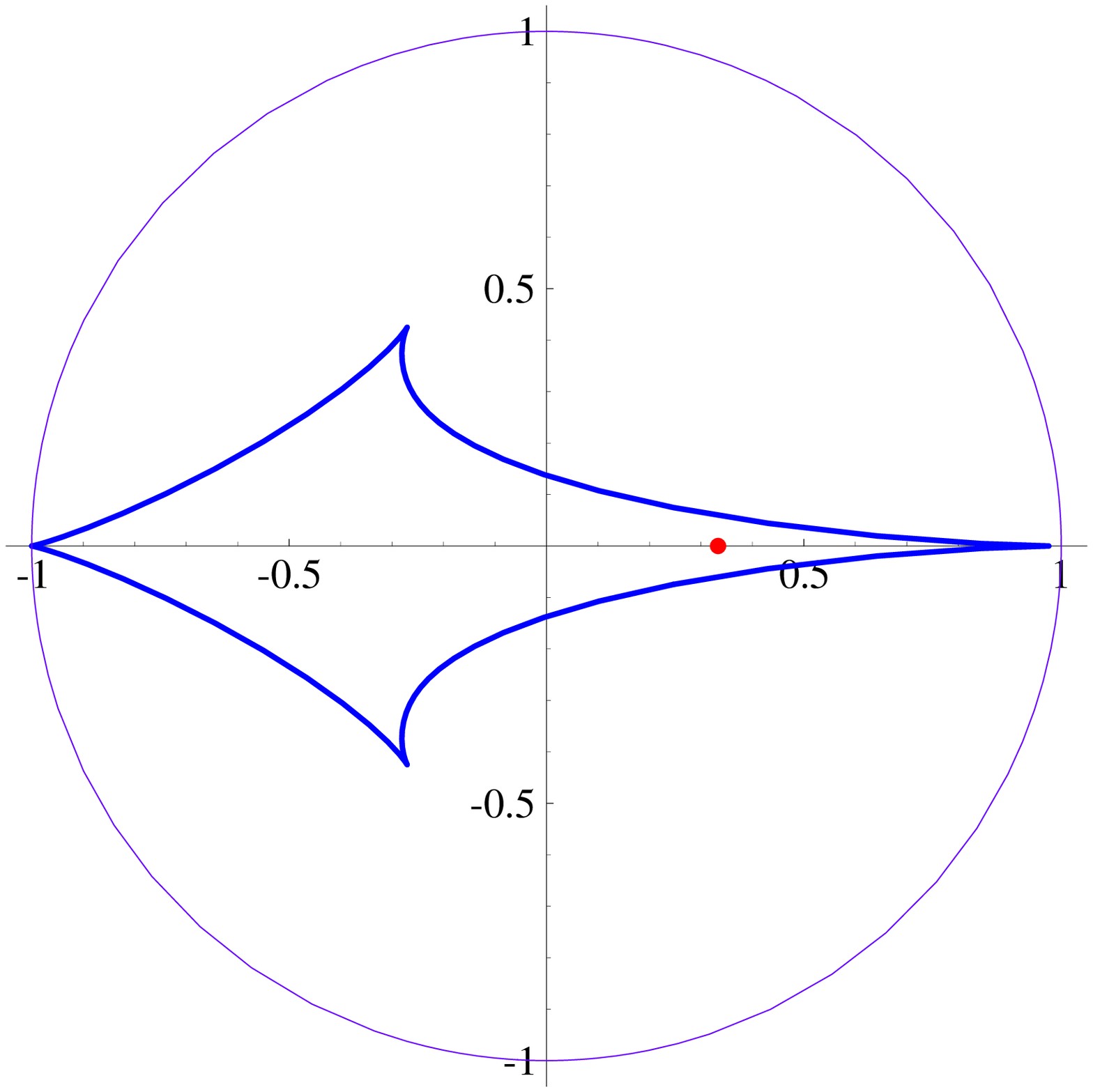}}\\
{\tiny caustic for $n=3$, $r=1/3$.}
} \hfill
\parbox{52mm}{\centering%
\mbox{\epsfxsize=50mm \epsfbox{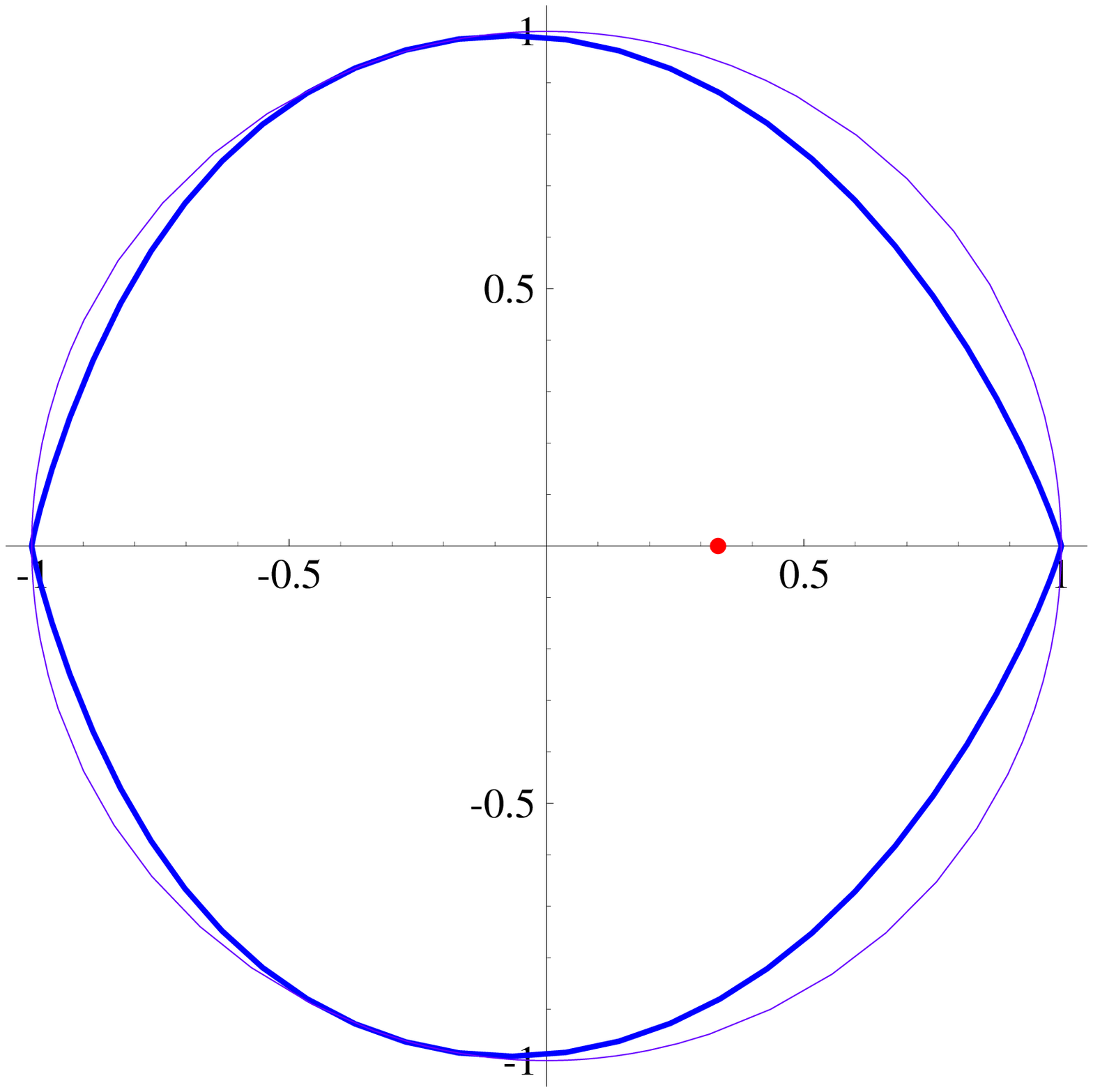}}\\
{\tiny caustic for $n=4$, $r=1/3$.}
}}
\medskip

\begin{proposition}\label{PropTgtReven}
Any caustic curve of $R_r^{2m}$ is tangential to the unit
circle at
$(\cos\phi,\sin\phi)$ at any point $\phi$ satisfying
$R_r^{2m}(\phi) = \phi$. In particular, this includes
the points $0$, $\pi$, and $\pm\phi_c$.
Moreover, if $r\leq 1/3$, these are the only four points
that the caustic
meets the circle.
\end{proposition}
\begin{proof}
By substitution of such $\phi$ into equations
(\ref{Eqncausticf}) with
$f=R_r^n$,
we have $x(\phi) = \cos\phi$ and
$y(\phi) = \sin\phi$.  Applying lemma~\ref{LemCuspAtFixes} to
such $\phi$'s, the assertion about the tangential property
of the caustics of even
iterations of $R_r$ follows easily.
Secondly, since
$x^2 + y^2 = 1$, we have
\begin{align*}
(R_r^n)'(\phi)\cos\phi - \cos R_r^n(\phi) &= 0 \\
(R_r^n)'(\phi)\sin\phi - \sin R_r^n(\phi) &= 0.
\end{align*}
This leads to $R_r^n(\phi) = \phi$.  If $r\leq 1/3$, by
lemma~\ref{LemCalRn}, $n$ must be even and $\phi$ is one 
of the four 2-periodic points.
\end{proof}
We have already seen from symplectic topology that caustics
of odd
iterations, $R_r^{2m+1}$, always have at least four cusps for
sufficiently small $r>0$.  Moreover,
for all $R_r$, two of the cusps occur at
$\phi=0,\pi$.
Now, we may extend this result to odd iterations of $R_r$
with isolated exceptional values of $r$.
It turns out these exceptional values only occur in
$r>1/3$.
\begin{proposition}\label{PropCusp0piRodd}
For $0<r< 1/3$ and for generic $1/3 < r < 1$, the caustic of
$R_r^{2m+1}$ always has cusps at $0, \pi$.
\end{proposition}
\begin{proof}
Note that for both $\phi_a = 0, \pi$, one has
$R_r(\phi_a) = \phi_a + \pi \mod 2\pi$, thus we may apply
lemma \ref{LemCuspAtFixes} to check whether there are cusps.
Furthermore, $R_r^2(\phi_a)=\phi_a$ and $R_r''(\phi_a)=0$,
we can simplify the chain rules (lemma~\ref{LemChainRuleRn}) and
obtain $(R_r^{2m+1})''(\phi_a)=0$ and
\begin{align*}
(R_r^{n+2})^{(3)}(\phi_a) &=
R_r'(R_r(\phi_a)) (R_r^n)'(\phi_a) R_r^{(3)}(\phi_a) + \\
&\qquad {} +
   {{R_r'(\phi_a)}^3} \left[ (R_r^n)'(\phi_a)
   R_r^{(3)}(R_r(\phi_a)) +
	 {{R_r'(R_r(\phi_a))}^3} (R_r^n)^{(3)}(\phi_a)
	 \right].
\end{align*}
Let us temporarily define, for positive integers $n$,
$$
A_{n} = -(R_r^{n})' + {(R_r^{n})'}^3 + 2 (R_r^{n})^{(3)}.
$$
So it suffices to check that $A_{2m+1}(\phi_a) \ne 0$.
We will proceed by induction on $m$.  First, by
direct computation,
$$
A_1(0) = \frac{24r^2}{(1-r)^2} > 0, \qquad\qquad
A_1(\pi) = \frac{24r^2}{(1+r)^2} > 0.
$$
Then, it can be shown that 
\begin{align*}
A_{n+2}(\phi_a) 
&= \left(\frac{1-9r^2}{1-r^2}\right)^3 A_n(\phi_a) +
\left[2R_r^{(3)}(\phi_a+\pi)R_r'(\phi_a)^3 +
2R_r'(\phi_a+\pi)R_r^{(3)}(\phi_a) \right. \\
&\qquad \qquad \left. {} - R_r'(\phi_a)R_r'(\phi_a+\pi) +
R_r'(\phi_a)^3R_r'(\phi_a+\pi)^3 \right]
(R_r^n)'(\phi_a) \\
&= \left(\frac{1-9r^2}{1-r^2}\right)^3 A_n(\phi_a) +
A_2(\phi_a)\cdot(R_r^n)'(\phi_a). \\
\intertext{In particular,}
A_{2m+1}(0) &= \left(\frac{1-9r^2}{1-r^2}\right)^3
A_{2m-1}(0) + A_2(0) \left(\frac{1-3r}{1-r}\right)^{m}
\left(\frac{1+3r}{1+r}\right)^{m-1}; \\
A_{2m+1}(\pi) &= \left(\frac{1-9r^2}{1-r^2}\right)^3
A_{2m-1}(\pi) + A_2(\pi) \left(\frac{1-3r}{1-r}\right)^{m-1}
\left(\frac{1+3r}{1+r}\right)^{m}; \\
\intertext{where}
A_2(0) &= \frac{48r^2(1+r)}{(1-r^2)^3}(1-3r+13r^2-15r^3), \\
A_2(\pi) &= \frac{48r^2(1+r)}{(1-r^2)^3}(1+3r+13r^2+15r^3).
\end{align*}
It can be easily computed that $A_2(0) > 0$ for $0<r<1/3$ and
$A_2(\pi)>0$ for all $r>0$.  Thus,
$A_{2m+1}(\phi_a) \geq \left(\dfrac{1-9r^2}{1-r^2}\right)^3
 A_{2m-1}(\phi_a) > 0$.

For $r > 1/3$, $A_{2m+1}$ may have zeros.  We can only
conclude from above that $A_{2m+1}$ is a rational function
in $r$ with denominator being a power of $(1-r^2)$.
Therefore, it has only isolated zeros and cusps at $0, \pi$ occur
for $r > 1/3$ generically.
\end{proof}
We observe that at cusp points, (\ref{Eqncausticftgt})
gives a set of ``homogeneous'' equations which has zero determinant.  
Thus, except at a couple of $\phi$'s, it is sufficient to solve
only one equation, say, $x'(\phi)=0$.  It is likely that this equation
has exactly four solutions for $0<r<1/3$.  However, it is still hard
to solve explicitly especially for high iterations of $R_r$.
Here are two pictures of the functions
$(R_r^3)''/((R_r^3)'(1+(R_r^3)'))$ and $(\sin\phi + \sin
R_r^3(\phi))/(\cos\phi - \cos R_r^3(\phi))$ for $r = 0.1$ and
$r=0.33$.  Exactly four cusp solutions are demonstrated in each of
them.  We
do not have a proof of this graphical fact.  Perhaps it
may be proved by detailed
curve sketching argument and comparison of $R_r^{2m-1}$ and
$R_r^{2m+1}$ using the known properties of $R_r^2$ given in
lemma~\ref{LemCalRn}.

{\parindent=0mm
\parbox{72mm}{\centering%
\mbox{\epsfxsize=70mm \epsfbox{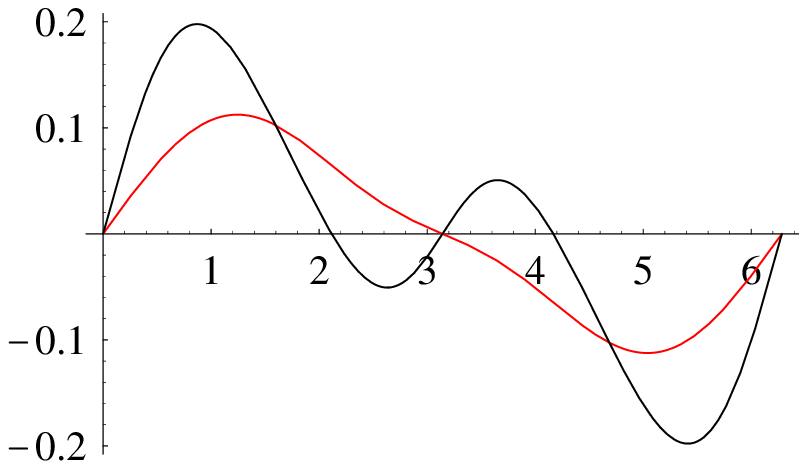}}
} \hfill
\parbox{72mm}{\centering%
\mbox{\epsfxsize=70mm \epsfbox{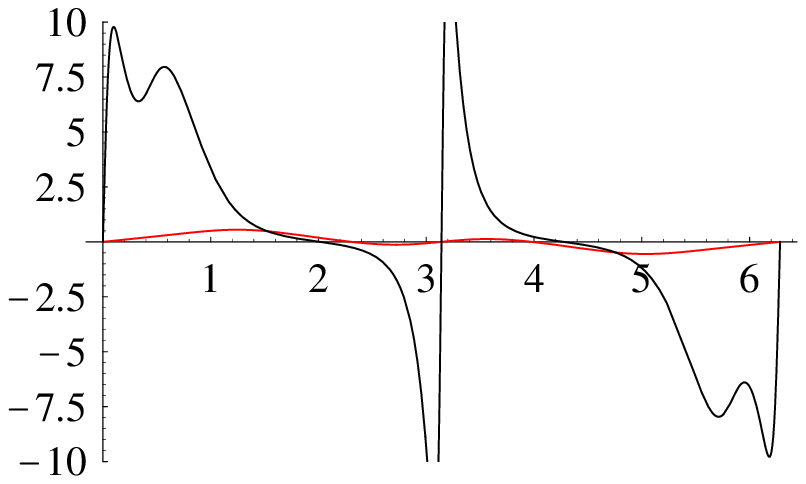}}
}
}

\noindent
It is probably worthwhile to compare the situation here with the
so-called Jacobi conjecture promoted by Arnold \cite{A1,A5}. As mentioned
in the beginning of \S\ref{Sec-caustics}, the caustics of 
$R_r$ agrees with the locus
of conjugate points of $(r,0)$ on a flat flying disk. Although
the loci of higher order conjugate points are not the same as the
caustics of odd iterations of $R_r$, their common contact geometric
nature
indicates that two problems of whether there are exactly 4
cusps on the loci of higher order conjugate points and whether
there are exactly 4 cusps of the caustics of odd iteration of $R_r$
might
be related. In both cases, since we do not have exact nice formulae
for the
loci of higher order conjugate points and the caustics of odd
iterations
of $R_r$, it would be very difficult to have an exact count of
cusps. Of course, Conjecture~\ref{ConjAllRr} deals with a 
very special situation. From graphical
evidence, it is tempting to think that a proof should not be out of
reach
with brute force calculation. Nevertheless, after spending much
effort on
this temptation, we think some more conceptual understanding of the
caustics of iterations of $R_r$ is needed in order to get an exact
count of cusps.

Furthermore, in the proof above, we see that for 
$n\geq 3$, except $A_3(\pi)$, the quantity
$A_{n}(0) = 0 = A_n(\pi)$ always at $r=1/3$.  There is a
possible structural change on the caustic of $R_r^{n}$
occurring at $r=1/3$ for $n\geq 3$.  
We observe from experiment that, for the caustic of odd 
iterations, once $r > 1/3$, bifurcation
of cusp may occur.  Interestingly, from the computed pictures (only
that of $R_r^3$ is shown), bifurcation only occur at the cusp
corresponding to $\phi=0$ but not others.  Would
the different properties between $A_2(0)$ and
$A_2(\pi)$ be part of the reasons?  

\medskip
{\parindent=0mm
\parbox{72mm}{\centering
\mbox{\epsfysize=50mm \epsfbox{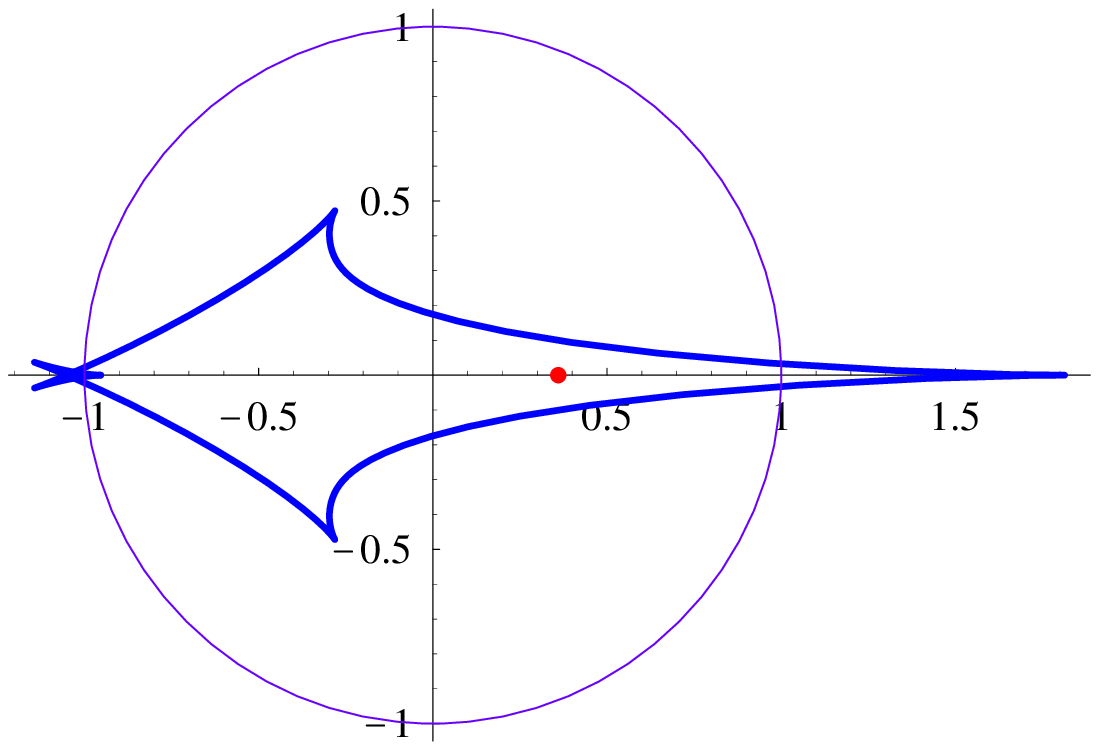}} \\
{\tiny $n = 3$, $r=0.36$.}
} \hfill
\parbox{80mm}{\centering
\mbox{\epsfysize=50mm \epsfbox{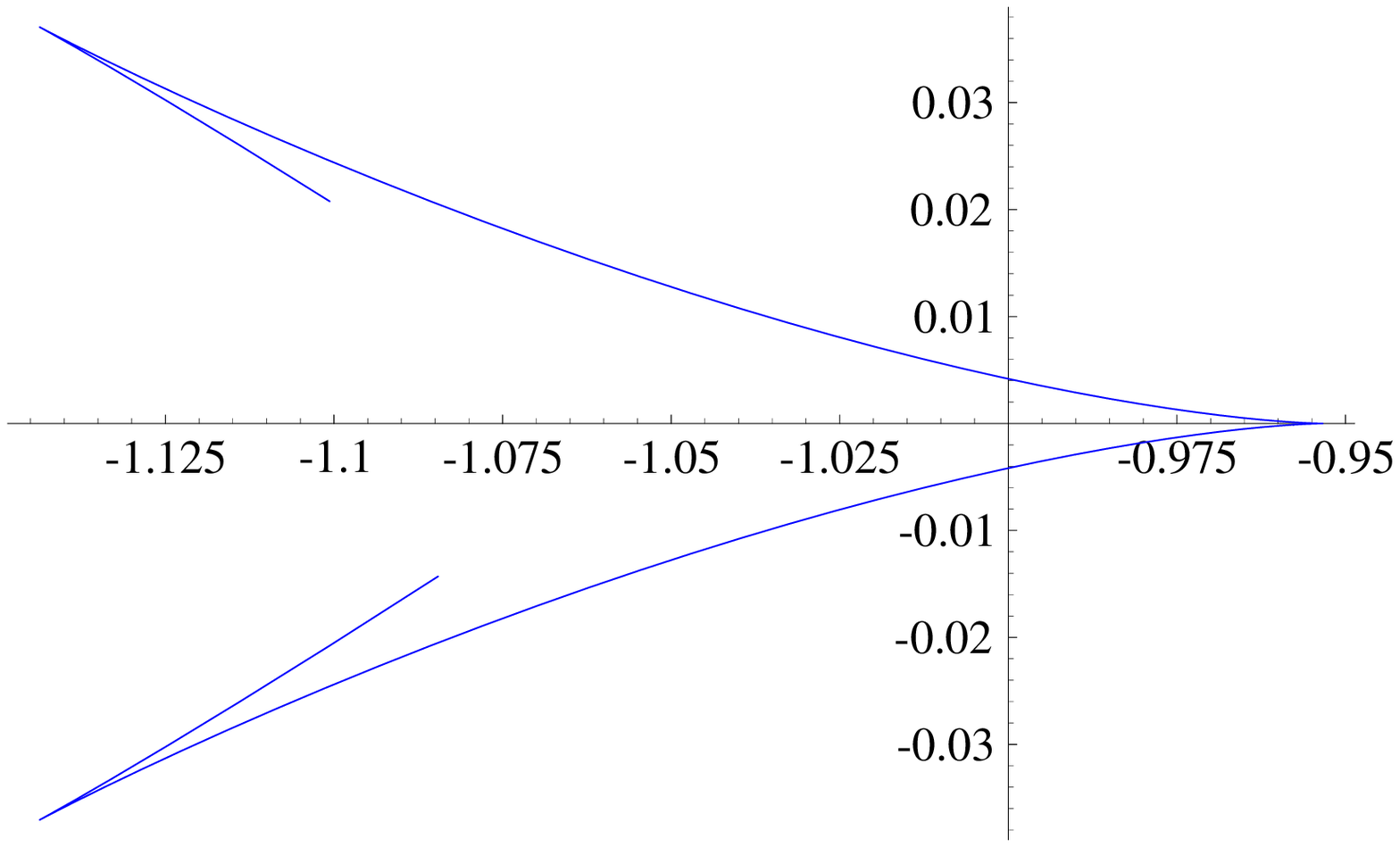}} \\
{\tiny its enlargement, $\phi=0$.}
}}
\medskip

On the other hand, bifurcation into cusps also occurs
for even iterations of $R_r$ at $r=1/3$.  
We will discuss this by beginning with some Taylor expansions.
Since $R_r$ and its iterations are $2\pi$-periodic odd functions, they
have particular nice expansions at $\phi_a=0,\pi$.  This enables us
to see the local properties of the caustics more clearly.

Let $f$ be any even iteration of $R_r$, then $f(\phi_a) = \phi_a$.
We write $\vartheta = \phi-\phi_a$ and $g(\vartheta) = f(\phi) -
\phi_a$ and suppose it has an expansion
$$
g(\vartheta) = \sum_{k=0} a_{2k+1}\vartheta^{2k+1}.
$$
One can inductively work out the coefficients of the expansions of
$1+f'(\phi) = 1+g'(\vartheta)$, $\cos f(\phi) = \pm\cos
g(\vartheta)$, etc.  If $P_k, Q_k$ denote
polynomials with $P_k(0,\ldots,0)=0=Q_k(0,\ldots,0)$, one has
\begin{align*}
x(\phi) &= \pm 1 \mp \frac{a_1}{2}\vartheta^2 + \sum_{k=1}
\frac{P_k(a_1,\ldots,a_{2k+1})}{(1+a_1)^{2k+1}} \vartheta^{2k}; \\
y(\phi) &= \pm \frac{2a_1}{1+a_1}\vartheta + \sum_{k=1}
\frac{Q_k(a_1,\ldots,a_{2k+1})}{(1+a_1)^{2k+1}}\vartheta^{2k+1}.
\end{align*}
These expansions are helpful to understand the caustics of
$R_r^{2m}$ at $\phi_a=0,\pi$.  It would be convenient to look at the
pictures before we go on.

\medskip
\noindent{\parindent=0mm
\parbox{52mm}{\centering
\mbox{\epsfysize=50mm \epsfbox{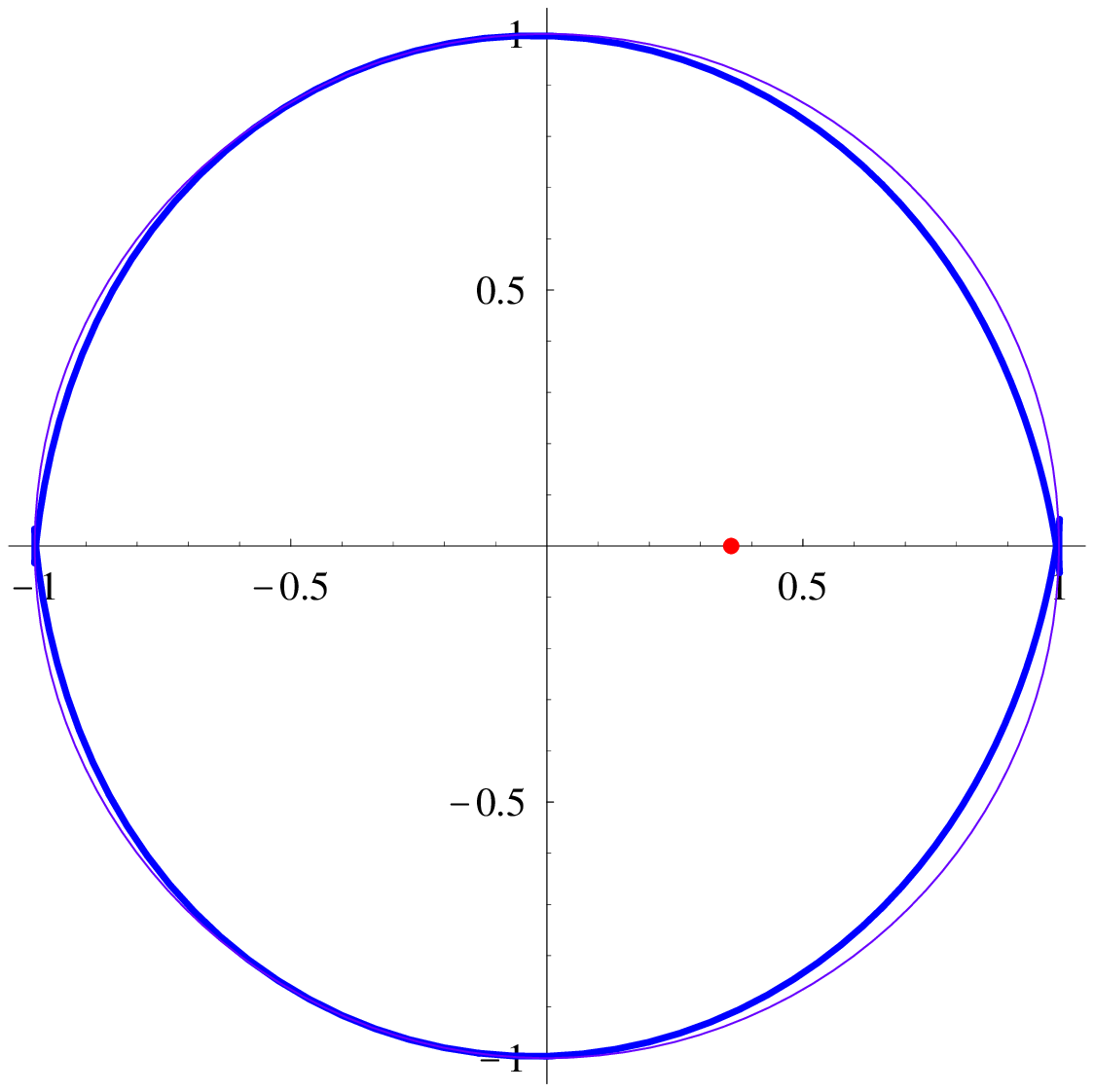}} \\
{\tiny $n = 2$, $r=0.36$.}
} \hfill
\parbox{73mm}{\centering
\mbox{\epsfysize=50mm \epsfbox{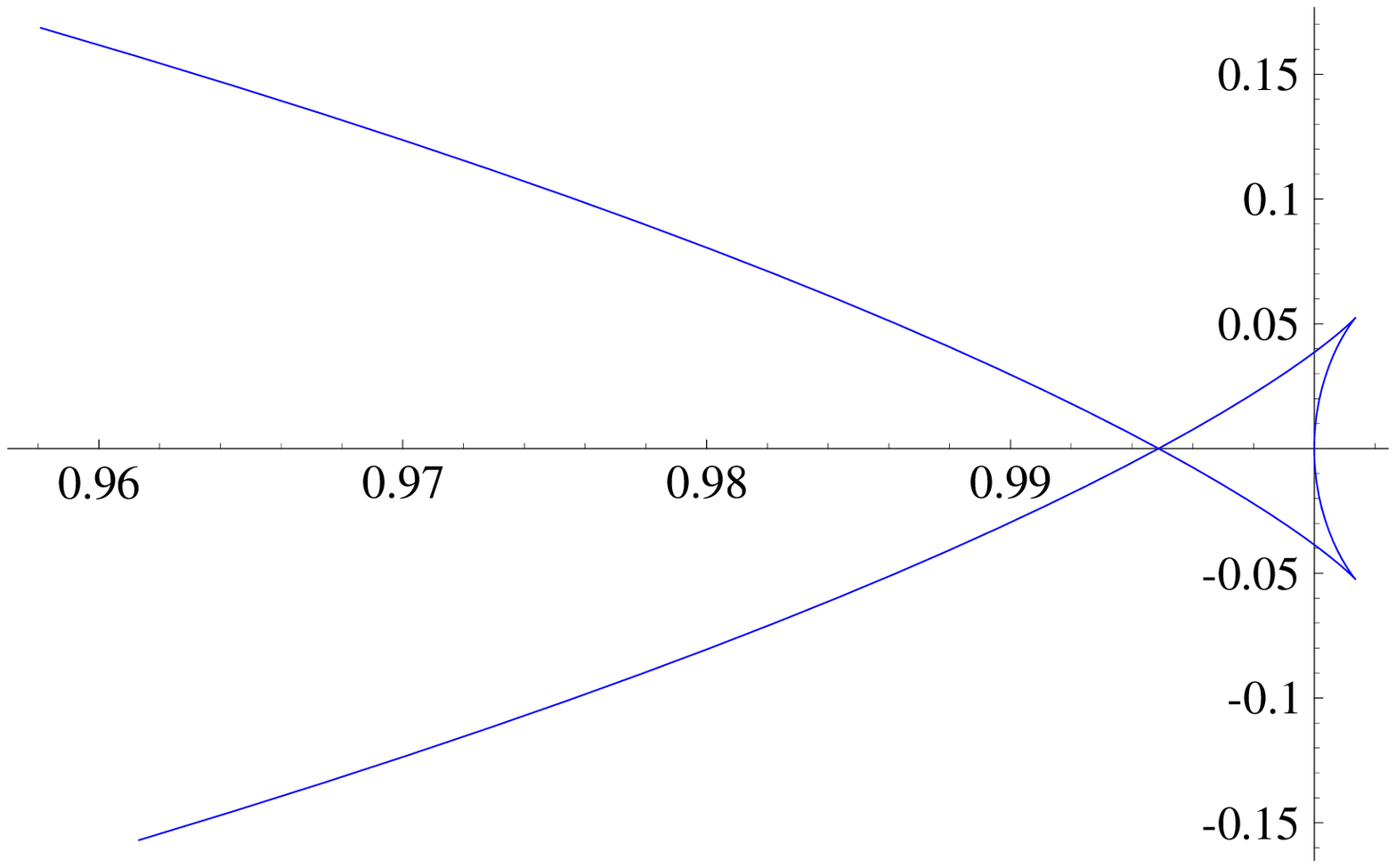}} \\
{\tiny its enlargement, $\phi=0$.}
}\hfill\\
\medskip
\noindent\parbox{52mm}{\centering
\mbox{\epsfysize=50mm \epsfbox{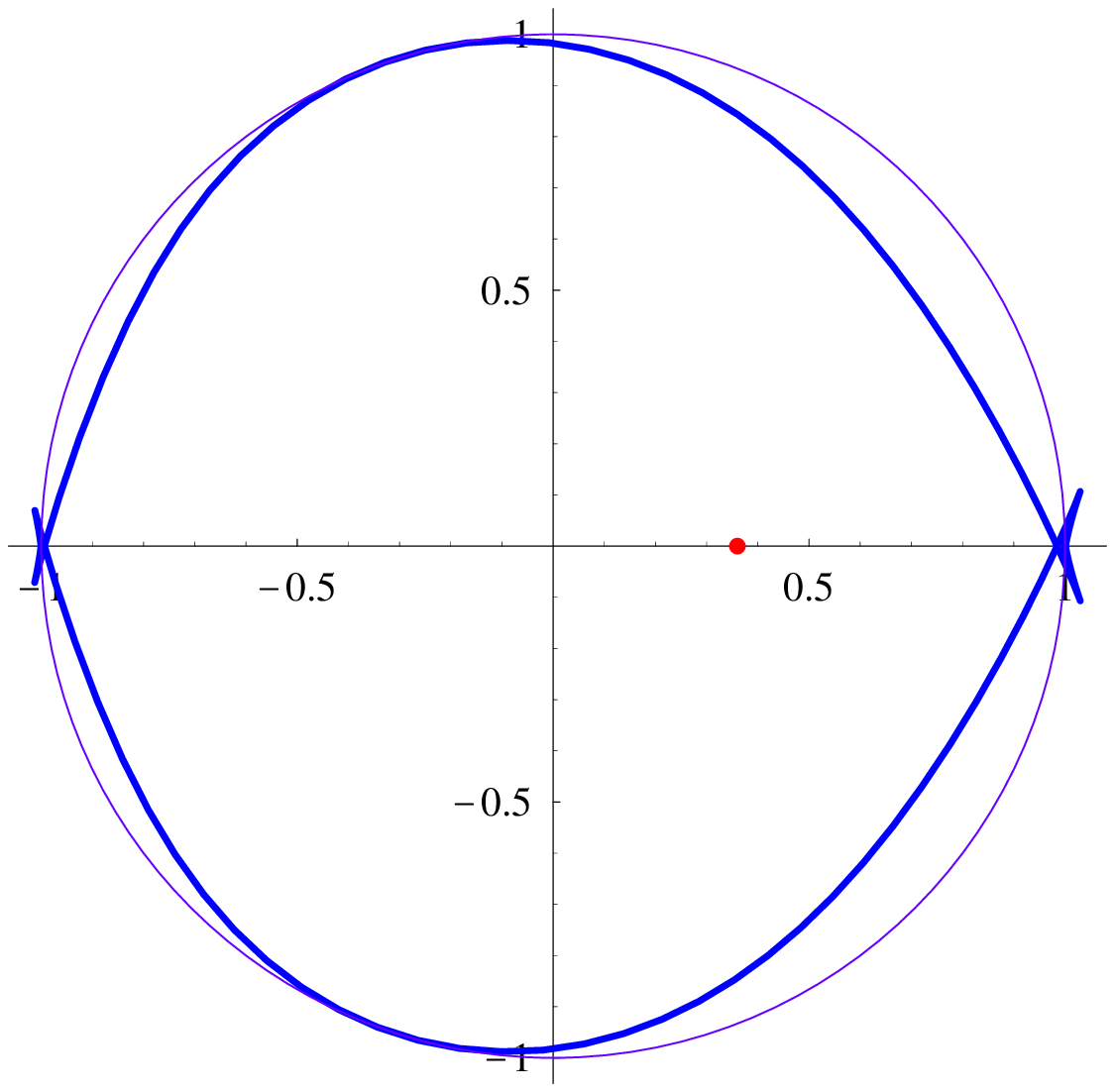}} \\
{\tiny $n = 4$, $r=0.36$.}
} \hfill
\parbox{73mm}{\centering
\mbox{\epsfysize=50mm \epsfbox{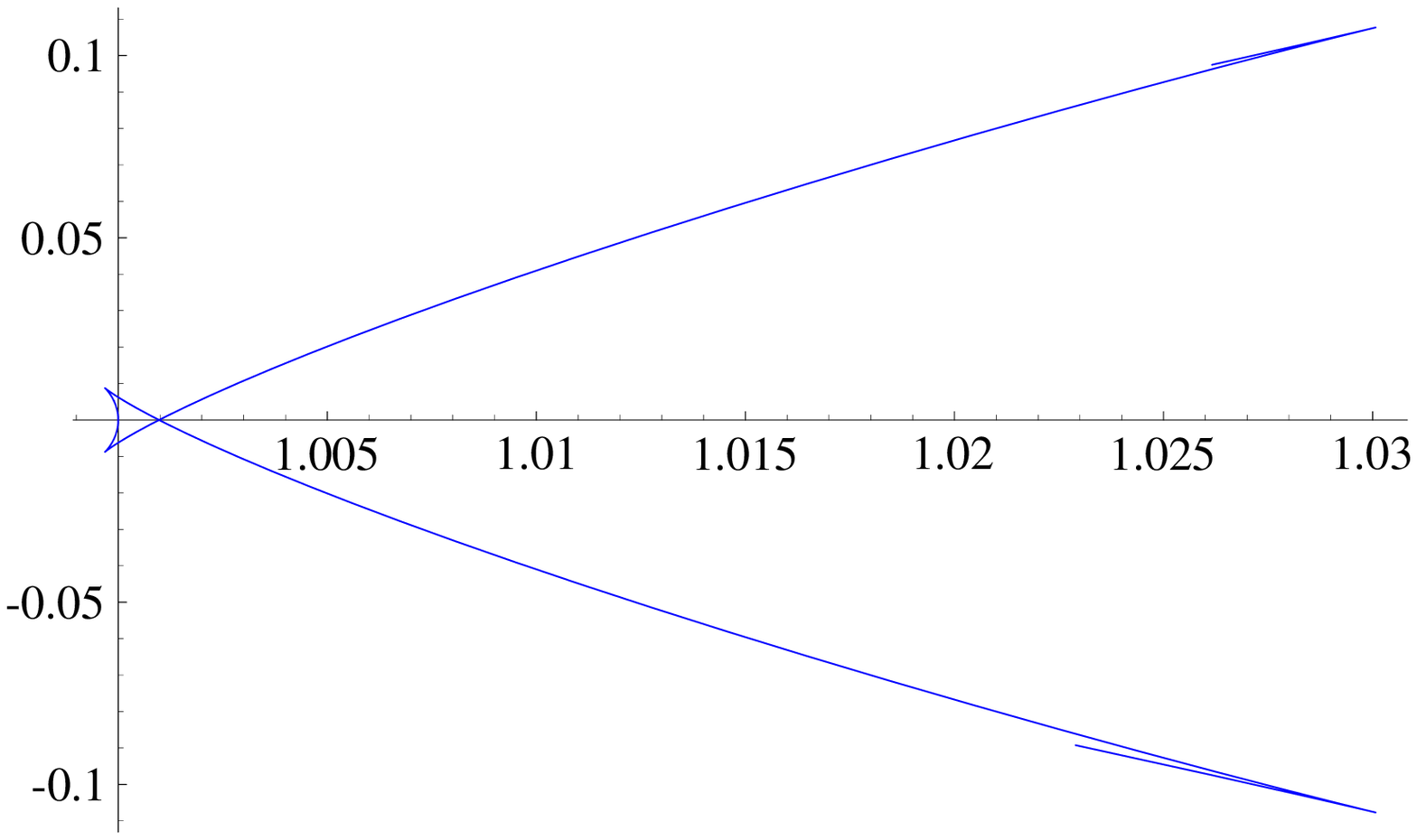}} \\
{\tiny its enlargement, $\phi=0$.}
}}
\medskip

\noindent In the above pictures, cusps are born near
$\phi=0$ and $\pi$.  From the enlargement, the
caustic bifurcates into $2m$~cusps when $r$ increases 
across $1/3$, where $R_r$ changes from a diffeomorphism to a
degree~1 map.  In the expansion of $R_r^{2m}$, $a_1 = \left(
\frac{1-9r^2}{1-r^2} \right)^m$.  By this, we will see that
$r=1/3$ is the value that the
caustic of $R_r^{2m}$ changes at $\phi=0,\pi$.
In fact, the caustics of
$R_{1/3}^2$ has the following Taylor expansions.  At $\phi_a = 0$,
\begin{align*}
x(\phi) - x(0) &=
\frac{-27}{4}\vartheta^4 + \text{O}(\vartheta^6) \\
y(\phi) - y(0) &=
18\vartheta^3 + \text{O}(\vartheta^5); \\
\intertext{and at $\phi_a=\pi$,}
x(\phi) - x(\pi) &=
\frac{243}{16}\vartheta^4 + \text{O}(\vartheta^6) \\
y(\phi) - y(\pi) &=
\frac{-81}{2}\vartheta^3 + \text{O}(\vartheta^5).
\end{align*}
This shows that the caustic of $R_r^2$ undergoes a
swallowtail bifurcation 
at $0,\pi$ when $r=1/3$.
We may further work out the
expansion of $R_{1/3}^{2m}$ as the
$m$-iteration of $R_{1/3}^2$. Using $a_1=0$ and $a_3\ne 0$
for $R_r^2$, we have
$$
R_{1/3}^{2m}(\phi) = \phi_a + \vartheta^{3^m} U(\vartheta)
$$
for some function $U$ with $U(0)\ne 0$.  The bifurcation of
the caustics of $R_r^{2m}$ at $\phi=0,\pi$
should be of the type $(\vartheta^{3^m+1},\vartheta^{3^m})$ when $r$
passes $1/3$.

In the above, we consider the behavior of the caustics $R_r^n$
with parameter $r$ and $n$ fixed.
What happens to the caustics if $r$ is 
fixed and $n$ is allowed to vary?  There is an interesting
phenomenon for $r \leq 1/3$.  Although
there is fundamental difference when $n$ is odd and even,  
this difference
disappears as $n$ goes to infinity.  
Here are the pictures showing how the caustics
of $R_r^{2m+1}$ and $R_r^{2m}$ change--they tend 
to the same quadrilateral.

\medskip
{\parindent=0mm
\parbox{52mm}{\centering
\mbox{\epsfxsize=50mm \epsfbox{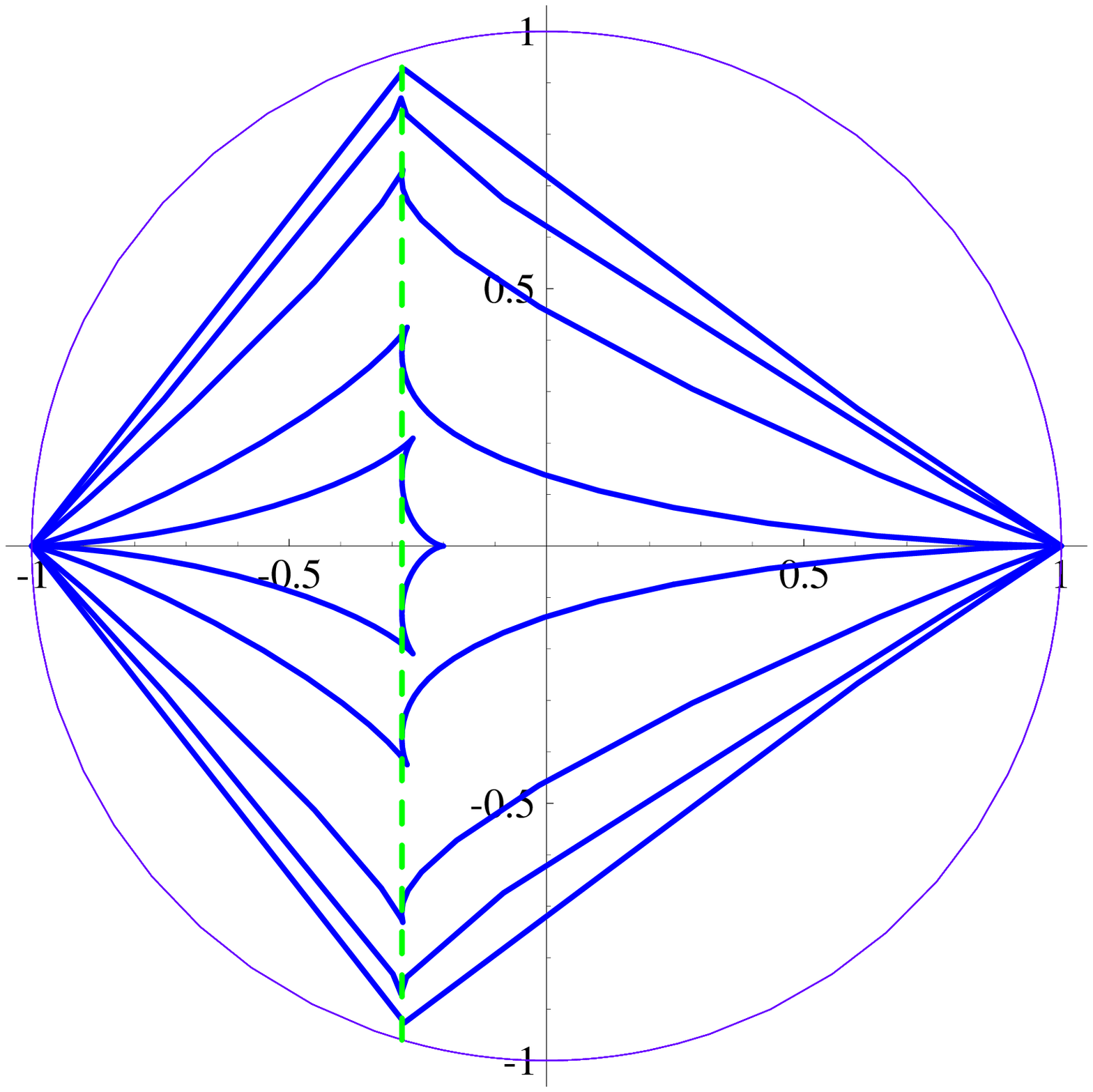}} \\
{\tiny caustics of $R_{1/3}^{\text{odd}}$}
}\hfill
\parbox{52mm}{\centering
\mbox{\epsfxsize=50mm \epsfbox{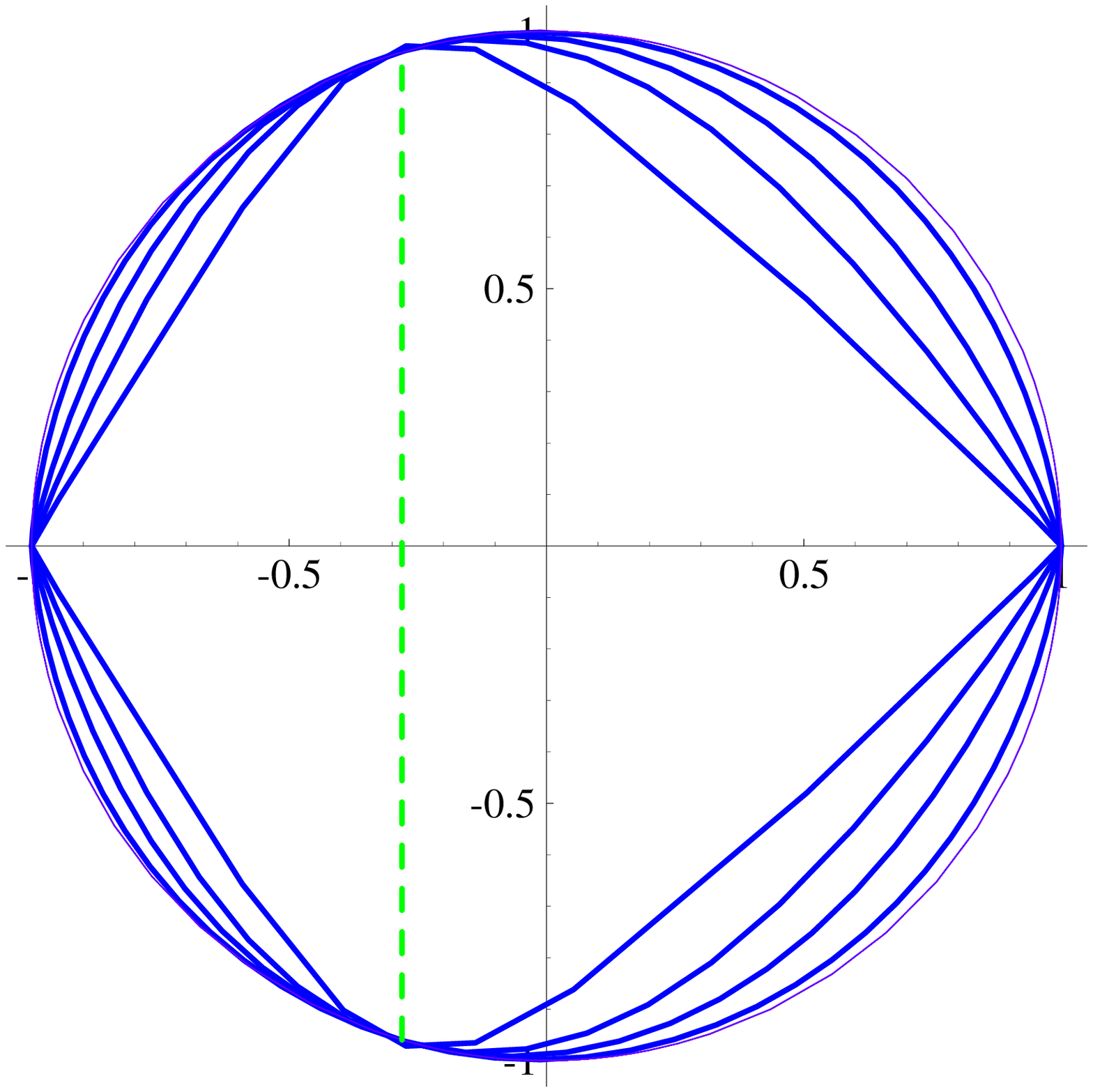}} \\
{\tiny caustics of $R_{1/3}^{\text{even}}$}
}}
\medskip

\noindent The green (dashed) vertical line from $\phi_c$ to
$-\phi_c$ in the picture of even or
odd caustics is determined by $r$ but not the number of
iterations.  When $n$ is even, it is where the caustic is
tangent to the circle.  When $n$ is odd, every caustic is
tangent to this vertical line because it is the line joining
$\phi_c$ and $R_r^{2m+1}(\phi_c) = -\phi_c$.  The point of tangency
occurs exactly at $\phi_c$ by definition.
\begin{proposition}\label{PropQuadRinfty}
For $0<r\leq 1/3$, as $m\to\infty$, both the caustics of 
$R_r^{2m+1}$ and $R_r^{2m}$ approach the same quadrilateral
defined by the four points $0$, $\pi$, and $\pm\phi_c$,
which are the only 2-periodic points of $R_r$.
\end{proposition}
\begin{proof}
Let us first see that the caustics of $R_r^{2m+1}$ at the four
points tend to the circle as $m\to\infty$.
These four $\phi$'s are the solution to $R_r^2(\phi) =
\phi$.  Thus,
we have
$(R_r^{2m+1})'(\phi) = R_r'(\phi)^{m+1}\cdot
R_r'(R_r(\phi))^m$.
 
At $\phi=0, \pi$, the coordinates of the caustic are
$$
x(\phi) =
\frac{\pm\left((R_r^n)'(\phi)-1\right)}
{1+(R_r^n)'(\phi)}, \qquad\qquad
y(\phi) = 0.
$$
It is clear that $(R_r^{2m+1})'(0) =
R_r'(0)^{m+1}R_r'(\pi)^m
= \left(\frac{1-3r}{1-r}\right)^{m+1}\cdot
\left(\frac{1+3r}{1+r}\right)^m$.  Thus, $x(0) \to -1$ as
$m\to\infty$.  The situation at $\pi$ is similar.
 
At the point $\phi_c$ with
$R_r^n(\phi_c) = -
\phi_c$, the coordinates of the caustic at this point are
$$
x(\phi_c) = \cos \phi_c, \qquad\qquad
y(\phi_c) = \frac{
\left(-1+(R_r^n)'(\phi_c)\right)\sin\phi_c}
{1+(R_r^n)'(\phi_c)}.
$$
Since $R_r$ is odd and $R_r^2(\phi_c) = \phi_c$, it follows
that
$(R_r^n)'(\phi_c) = R_r'(\phi_c)^n$.  Moreover, by
$4r\cos\phi_c =
1-\sqrt{1+8r^2}$,  one may show that
$$
\frac{\left(-1+(R_r^n)'(\phi_c)\right)}{1+(R_r^n)'(\phi_c)}
\to 1
\qquad\quad \text{as $n\to\infty$.}
$$
Again, these two cusps approach to the unit circle.

Secondly, from the lemmas~\ref{LemCalR} and~\ref{LemCalRn}, $0$ and $\pi$
are the attracting fixed points of $R_r^2$ while $\pm\phi_c$
are repelling.  Moreover, the attracting basins for $0$ and $\pi$ are
$(-\phi_c,\phi_c)$ and $(\phi_c,2\pi-\phi_c)$
respectively.  Thus, for any given neighborhood of $0$,
for sufficiently large $m$, for any neighboring $\phi_1, \phi_2 \in
(-\phi_c,\phi_c)$, $R_r^{2m}(\phi_1)$ and
$R_r^{2m}(\phi_2)$ lie in that neighborhood of $0$.
Hence, the intersection of the lines from $\phi_j$ to
$R_r^{2m}(\phi_j)$ lies in a neighborhood of the
quadrilateral.  The proof for the cases at $\pi$ and of odd
iterations are similar.
\end{proof}

%
\subsection{Technical Results}
\label{Sec-details}
\label{Sec-technical}
In this section, a couple of technical results will be
given.  They are mostly done by direct computations and the
methods may not be insightful.  However, they may be the
necessary evil for they will be
used to justify our observations in \S\ref{Sec-observations}.
The first one deals with
the existence of cusp at certain special
``symmetric'' positions.
\begin{lemma}\label{LemCuspAtFixes}
Let $f$ denote any iteration of $R_r$.
On the caustic of $f$, the conditions for the occurrence of 
a semicubical cusp at $\phi$ are
\begin{itemize}
\item
$f'(\phi) = 0$ and $f''(\phi) \ne 0$ if $f(\phi)=\phi$;
\item
$f''(\phi) = 0$ and
$-f'(\phi) + f'(\phi)^3 + 2f^{(3)}(\phi) \ne 0$ if $f(\phi) =
\phi+\pi$;
\end{itemize}
\end{lemma}
\begin{proof}
This is proved by computing the derivatives of 
(\ref{Eqncausticf}) and (\ref{Eqncausticftgt}), then evaluate at the
particular values $\phi_0$ or $\phi_a$.  
Then the result follows by verifying $x'=0=y'$ and
$x''y'''-x'''y'' \ne 0$.  
\end{proof}
{\em Remarks.\/}
Although this lemma is stated for an iteration of $R_r$, it is
actually true for any circle map.  Likewise, many results in
this section hold in a more general setting but we would 
like to focus on iterations of $R_r$.
Furthermore, at a point $\phi$ with
$f(\phi) = -\phi$, we always have $x'(\phi)=0$.  The
conditions are
\begin{gather*}
f'(\phi)(1+f'(\phi))\cos\phi + f''(\phi)\sin\phi = 0 \\
f'(\phi)(2+\cos(2\phi)) + 6f'(\phi)^2\cos^2\phi +
(1+2\cos(2\phi))f'(\phi)^3 - 2f^{(3)}(\phi)\sin^2\phi \ne 0.
\\
\intertext{Analogously,  if $f(\phi) = \pi-\phi$, we have
$y'(\phi)=0$ and conditions}
f'(\phi)(1+f'(\phi))\sin\phi - f''(\phi)\cos\phi = 0 \\
f'(\phi)(2-\cos(2\phi)) + 6f'(\phi)^2\sin^2\phi +
(1-2\cos(2\phi))f'(\phi)^3 - 2f^{(3)}(\phi)\cos^2\phi \ne 0.
\end{gather*}
The next one may be an exercise for calculus students.
\begin{lemma}\label{LemChainRuleRn}
The chain rules for $R_r^n$ with $n=p+q$ are given by
\begin{align*}
(R_r^n)'(\phi) &= (R_r^p)'(R_r^q(\phi)) \cdot (R_r^q)'(\phi)
= R_r'(\phi)\cdot R_r'(R_r(\phi)) \cdot \cdots \cdot
R_r'(R_r^{n-1}(\phi))
\\
(R_r^n)''(\phi) &=
(R_r^p)''(R_r^q(\phi)) \cdot (R_r^q)'(\phi)^2 +
(R_r^p)'(R_r^q(\phi)) \cdot (R_r^q)''(\phi) \\
(R_r^n)^{(3)}(\phi) &= (R_r^p)^{(3)}(R_r^q(\phi)) \cdot
(R_r^q)'(\phi)^3 +
(R_r^p)'(R_r^q(\phi)) \cdot
(R_r^q)^{(3)}(\phi) + \\
{}&\qquad\quad + 3 (R_r^p)''(R_r^q(\phi)) \cdot (R_r^q)'(\phi)
\cdot (R_r^q)''(\phi).
\end{align*}
\end{lemma}
In determining the cusps on the caustic, some orbits in the
iteration play a special role.  We thus establish the
following to handle that.
\begin{lemma}\label{LemCalR}
For $0\leq r \leq 1/3$, $R_r^{2m+1}$ has no fixed 
point and $R_r^{2m}(\phi) \ne \phi+\pi$.
\end{lemma}
\begin{proof}
Firstly, one can obtain algebraically the four fixed points of
$R_r^2$.  The attracting ones are $0$, $\pi$, while $\pm\phi_c$ are
repelling.  Then by simple calculus, we get the
corresponding attracting basins and conclude that
$R_r^{2m}(\phi)$ converges to $0$ or $\pi \mod 2\pi$ monotonically. 
By the series expression (\ref{EqnRrSeries}) of $R_r(\phi)$, 
one can deduce the following estimate
$$
\left| R_r(\phi) - \phi \right| \geq \pi - 2\left| \log (1-r) \right| .
$$
Then, according to the convergence of  $R_r^{2m}(\phi)$, 
the lemma follows.
\end{proof}
\begin{lemma}\label{LemCalRn}
For $0<r\leq 1/3$, let $\phi_c =
\arccos\left(\dfrac{1-\sqrt{1+8r^2}}{4r}\right)$, 
the solution sets to the following
equations are given as,
\begin{itemize}
\item
$R_r^n(\phi)=\phi$ has solutions if and only if $n$ 
is even, which are $0, \pi, \pm\phi_c$.
\item
$R_r^n(\phi) = -\phi$ has solutions 
$\pm\phi_c$ if $n$ is odd; and $0, \pi$ if $n$ is even.
\item
$R_r^n(\phi)=\phi+\pi$ has solutions if and only if $n$ is odd, two of them are
$0$ and $\pi$.
\item
$R_r^n(\phi) = \pi-\phi$ has solutions
$0, \pi$ when $n$ is odd and
two solutions when $n$ is even.
\end{itemize}
\end{lemma}
\begin{proof}
First of all, the preceding lemma already give us part of
the conclusions.
Furthermore, $R_r^2$ fixes the intervals $[-\pi, -\phi_c]$,
$[-\phi_c,0]$, $[0,\phi_c]$, and
$[\phi_c,\pi]$ in a way that
\begin{align*}
R_r^2(\phi) &< \phi \qquad \text{on $(-\pi, -\phi_c)$ or $(0,\phi_c)$,} \\
R_r^2(\phi) &> \phi \qquad \text{on $(-\phi_c,0)$ or $(\phi_c,\pi)$.}
\end{align*}
Thus, induction process like $R_r^{2m}(\phi) < R_r^{2m-2}(\phi) < \phi$ on the
corresponding intervals will give the only fixed points of $R_r^{2m}$.
The other claims are done similarly.
\end{proof}
The above lemmas are illustrated by the following figures.

\medskip
{\parindent=0mm
\parbox{78mm}{\centering
\mbox{\epsfxsize=77mm \epsfbox{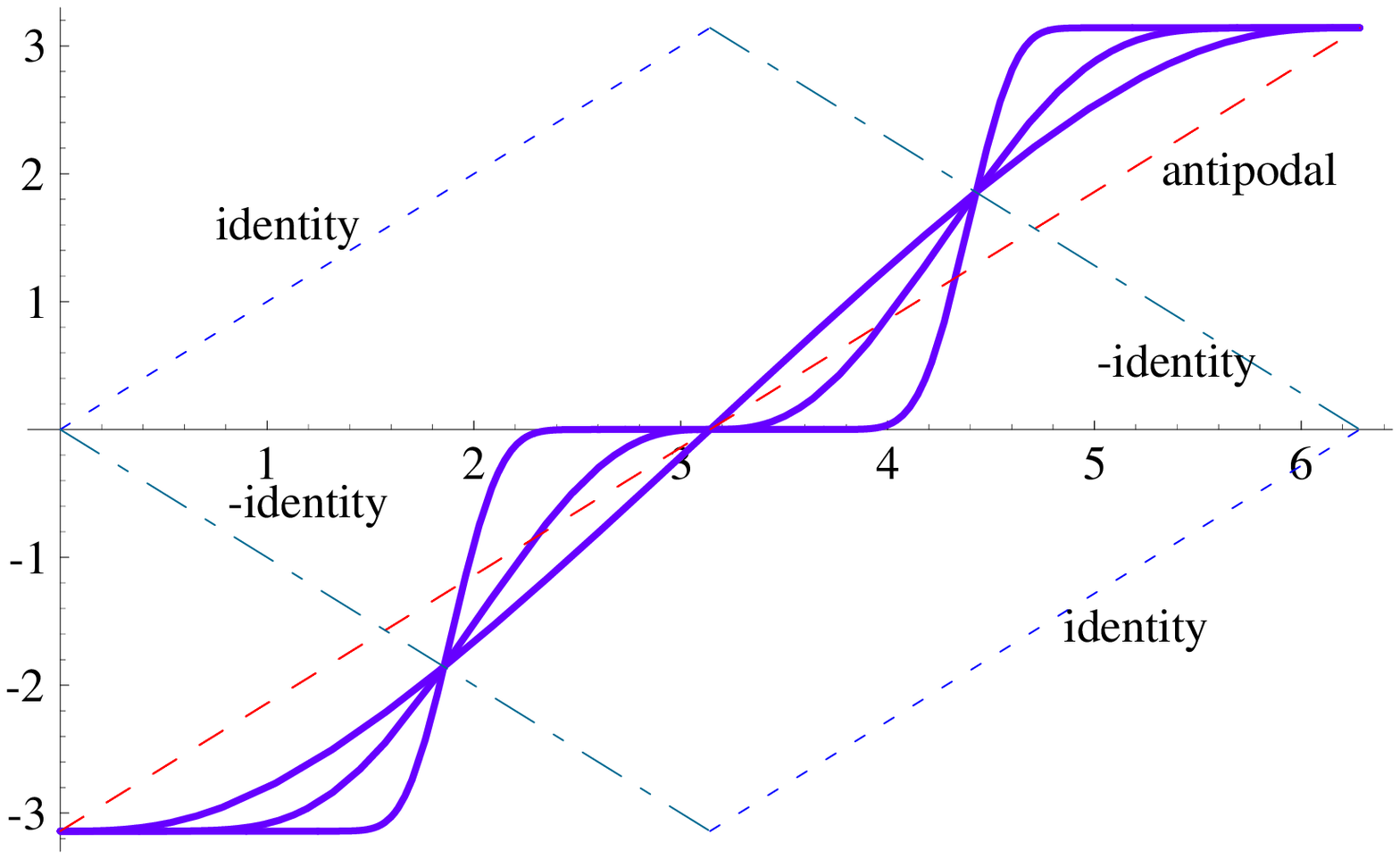}} \\
{\tiny Plots of  $R_r^1$, $R_r^3$ and $R_r^7$, $r=1/3$.}
}\hfill
\parbox{78mm}{\centering
\mbox{\epsfxsize=77mm \epsfbox{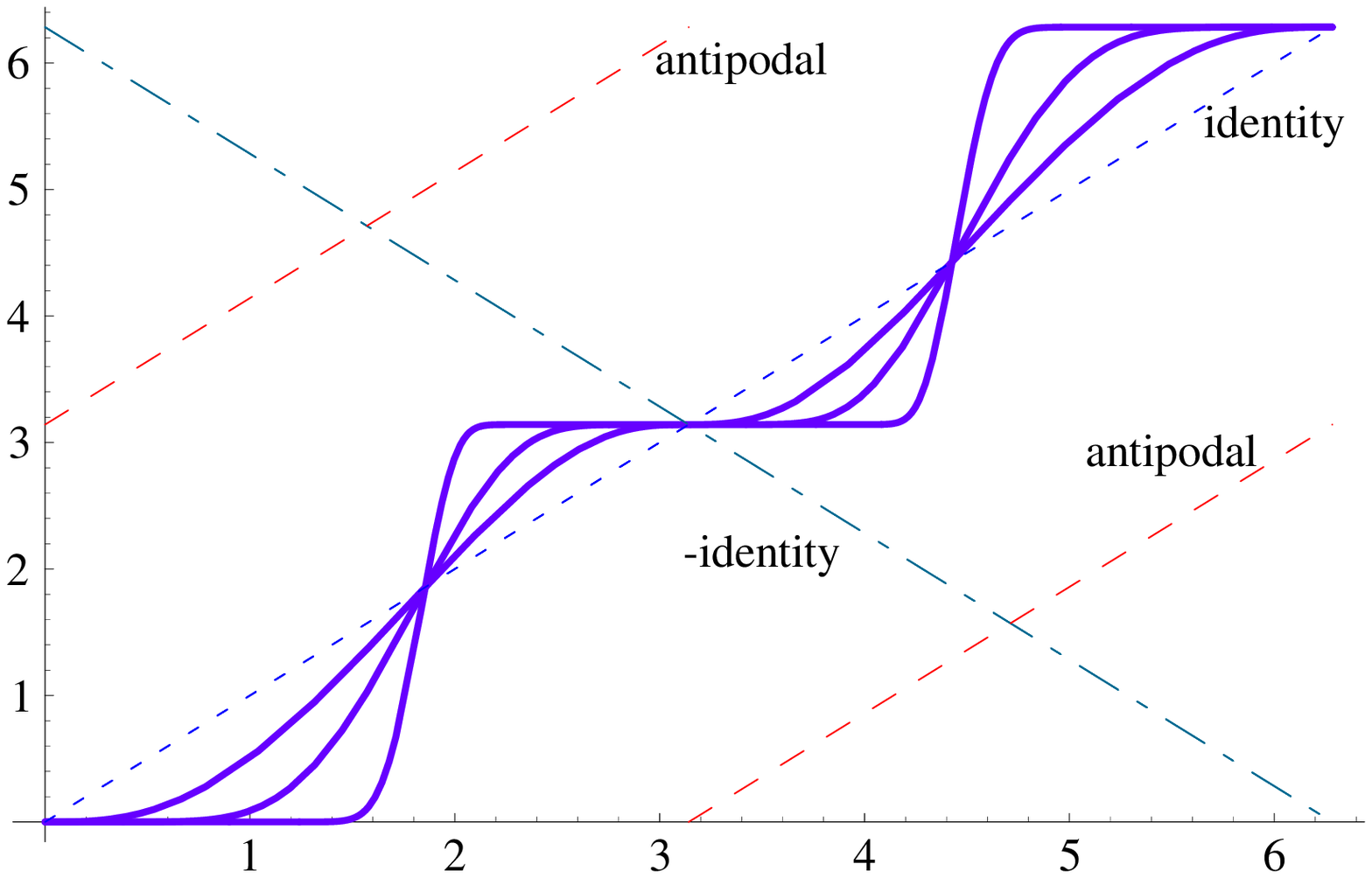}} \\
{\tiny Plots of $R_r^2$, $R_r^6$ and $R_r^8$, $r=1/3$.}
}}
\medskip

As an attempt to understand more about the iterations $R_r^n$, we
computed the {\sl asymptotic orbits} of the two critical points 
$\pm\arccos\left(\frac{1+3r^2}{4r}\right)$ for
$1/3 \leq r < 1$.
The asymptotic orbit of $\phi$ is the set
$\Set{R_r^n(\phi)~:~ N_1 < n < N_2}$ for large $N_1, N_2$.
The plot of asymptotic orbits against the parameter $r$ is called a
{\sl bifurcation diagram}.
It shows the attracting periodic cycles or chaotic behavior of a
map according to the variation of the parameter.  These diagrams
have been further analyzed in \cite{Au}.

\noindent
\parbox{77mm}{\centering
\mbox{\epsfxsize=75mm \epsfbox{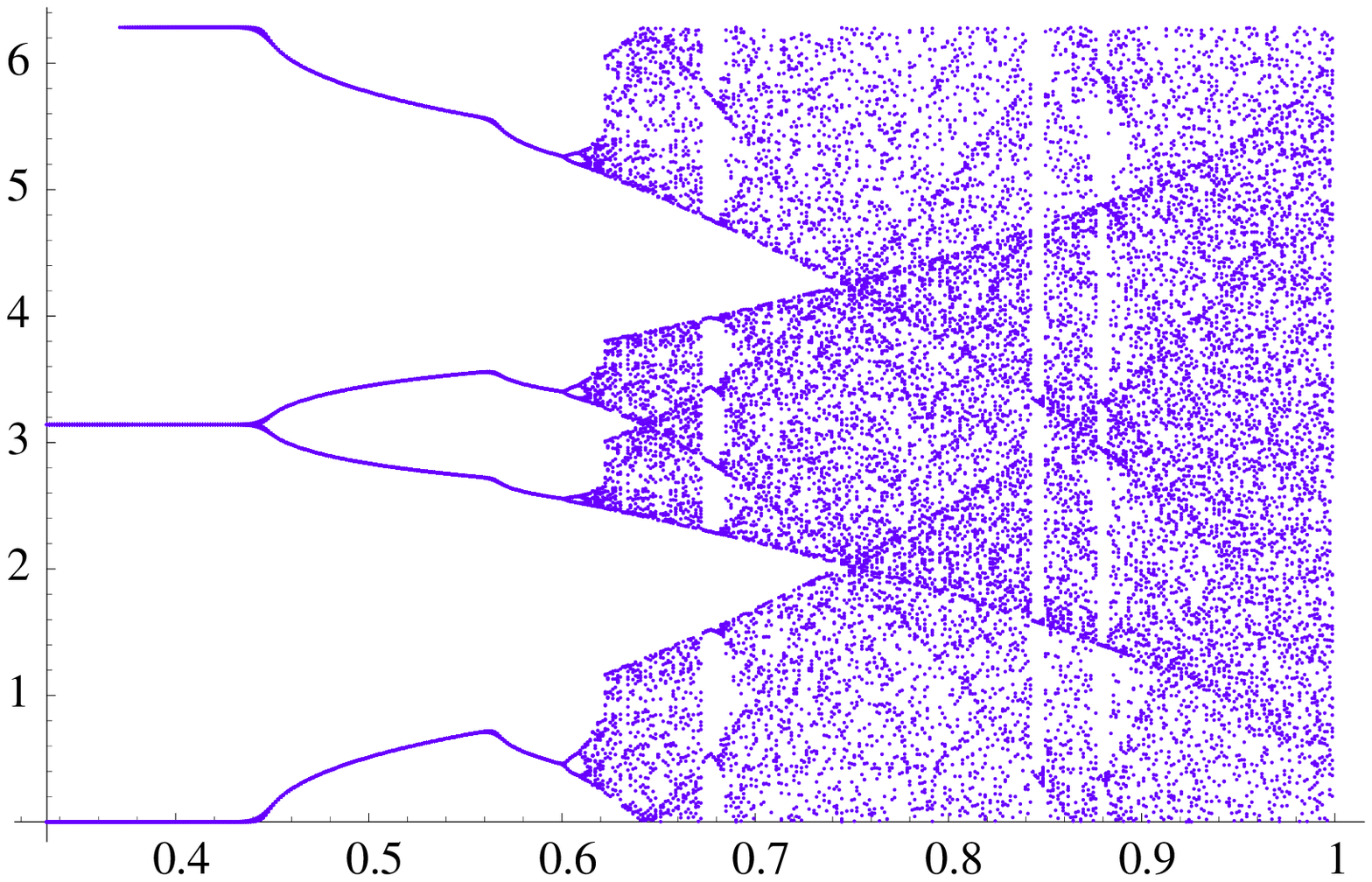}} \\
{\tiny asymptotic orbits of $\arccos\frac{1+3r^2}{4r}$}
}
\parbox{77mm}{\centering
\mbox{\epsfxsize=75mm \epsfbox{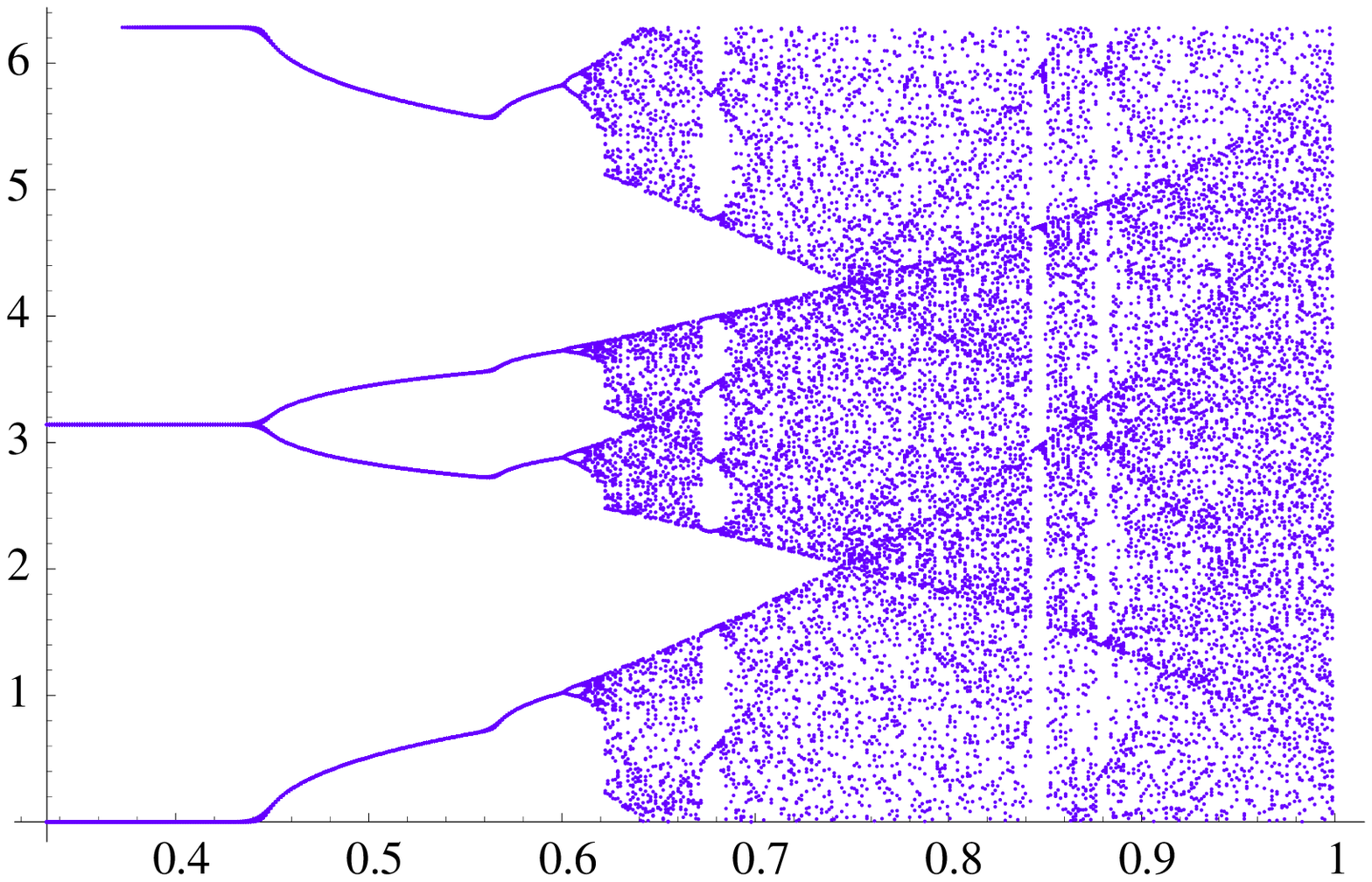}} \\
{\tiny asymptotic orbits of $-\arccos\frac{1+3r^2}{4r}$}
}

\section{Mode-locking}
\label{Sec-modelock}
\subsection{Background}
\label{Sec-background}
The study of circle maps is closely related to the study of differential equations on
torus (i.e., equations with double periodic coefficients).  For any such differential
equation, one may consider the Poincar\'{e} return map of the flow, which defines
a map on a meridian circle of the torus.  It turns out that the stability of the
equation is reflected by this circle map.

Far back in 1959, in his Ph.D. project, Arnold investigated the circle map
$$
\phi \mapsto \phi + a + \varepsilon \cos\phi
$$
and obtained information on its resonance zone in the
$(a,\varepsilon)$-plane, \cite[\S 12]{A2}.  This gives rise to the famous picture
of so-called Arnold tongues.  Subsequently, there are numerous studies, by
physicists and mathematicians, \cite{BBJ,Di,FKP,JBB,K,P,Z}, on the perturbation of a rotation
$$
\phi \mapsto \phi + \Omega - \varepsilon\sin\phi, \qquad \varepsilon\in[0,1).
$$
The focus is on the phenomenon called mode-locking and the
Devil's staircase.
Arnold later gave a proof in \cite{A4} of his observation
for circle maps of the form
$$
\phi \mapsto \phi + \Omega + \varepsilon (
\text{trigonometric polynomial})
$$
as well as analytic reduction of many circle maps, \cite[Ch.~3, \S
12]{A3}.
The algebraic nature of the method is also
apparent in the problem of particular differential equations.  
Arnold predicts that a
general theorem exists for these equations and general circle maps.

In this section, we will provide further evidence towards Arnold's
prediction by showing similar
behavior in the off-center reflection.  
It should be remarked that the off-center
reflection is not of the form studied by Arnold.
Thus, it may be another small step towards the general theory.

We consider a two-parameter model of circle maps which
arises from the off-center reflection map, namely, with parameters
$r\in [0,1)$ and $\Omega \in (-\pi,\pi]$,
$$
R_{r,\Omega}(\phi) = \phi + \Omega - 2\sum_{k=1}^\infty
\frac{r^k}{k}\sin(k\phi).
$$
Here we use $r$ instead of $\varepsilon$ to be consistent with 
previous sections. 
Note that unlike the models discussed above, $r$ 
cannot be factored out.
This map can be thought of as an imperfect off-center reflection on the
circle where the reflected angle has a constant deviation from the incident angle.  It
is the original off-center reflection when $\Omega=\pi$.  
We may not get such a deviation by varying the
metric of the circle; 
it is better understood in terms of symplectic geometry.

For $\phi_0\in\sphere^1$, there is the {\sl rotation number}
$\displaystyle
\omega(R_{r,\Omega},\phi_0) = \lim_{n\to\infty}
\frac{R^n_{r,\Omega}(\phi_0)-\phi_0}{n}$, where the
right hand side is performed on a lifting of
$R_{r,\Omega}$.  It is independent of
$\phi_0$ if $R_{r,\Omega}$ is diffeomorphic.  In such
case, one simply denotes $\omega(R_{r,\Omega})$.  If
$R_{r,\Omega}$ is only a degree~1 map, one has a {\sl
rotation interval} instead.  These notions are indeed defined for
any circle map.  Historically, attention has been centred around
perturbations of rotations,
$\phi \mapsto \phi + \Omega + u(\phi)$.
It is natural to ask for the relation between $\Omega$ and $\omega$.
The physicists usually refer to $\Omega$ as internal frequency and
$\omega$ as resonance frequency.  When $\omega=\omega(\Omega)$ is a
locally constant function, the situation
is called {\sl mode-locking}.

Herman has extensively studied the mode-locking property and
obtained interesting results, \cite{H1,H2}.  These
results are applicable to $R_{r,\Omega}$ because it
satisfies the property ${\bold A}_0$ defined by Herman.
\begin{theorem}\label{ThmOCmodelock}
For all $\omega_0 \in 2\pi\Q$ and $0< r\leq 1/3$, there is an interval
${\cal I} = {\cal I}_r$ of $\omega_0$ such that for every
$\Omega\in{\cal I}_r$, the diffeomorphism $R_{r,\Omega}$
has rotation number
$\omega_0$.
\end{theorem}
The interval ${\cal I}$ is called {\sl resonance interval} and
its size depends on $r$ (and of course
$\omega_0$).  
Its variance in terms of
$r$ defines a picture which looks like a tongue.  
We will discuss it later.
Furthermore, from Herman's study, the off-center 
reflection model also demonstrates the well-known
Devil's staircase.
\begin{theorem}\label{ThmOCstairs}
For any $0<r\leq 1/3$, the function $\Omega \mapsto
\omega(R_{r,\Omega})$ is nondecreasing, locally constant at any
rational number, and has a Cantor set of discontinuity.
\end{theorem}


We have mentioned that if we alter a sign and form the ``conjugate'' family
$$
\overline{R}_{r,\Omega}(\phi) =
\phi + \Omega + 2\sum_{k=1}^\infty \frac{r^k}{k}\sin(k\phi),
$$
the dynamics is completely different.  Actually,
$\overline{R}_{r,\Omega}$ can be extended to
$e^{2\pi\thei\Omega}\dfrac{z-r}{1-rz}$ on the hyperbolic disk, which defines a
hyperbolic element in $\text{PSL}(2,\R)$.  The mode-locking phenomenon does
not occur, i.e., $\omega\left(\overline{R}_{r,\Omega}\right)=2p\pi/q$
only if $\Omega = p/q$.

\subsection{Width of Resonance Zone}
\label{Sec-width}
In \cite{A4}, Arnold discusses the mode-locking situation of a rotation
slightly perturbed by a trigonometric polynomial, $g(x)$,
$$
f~:~x \mapsto x + \Omega + \varepsilon g(x).
$$
The {\sl resonance zone} is the set
$\Set{(\Omega,\varepsilon)~:~ \Omega \in {\cal I}_\varepsilon}$.
Arnold developed a formal calculation to estimate the width of the
interval ${\cal I}_\varepsilon$ in terms of $\varepsilon$, which gives
rises to a
picture of the resonance zone.  This formal calculation is related to the
homological equation of analytical reduction, \cite{A3}.
If the rotation number is rational, the width of the 
resonance interval ${\cal I}_\varepsilon$ is bounded by a
power of $\varepsilon$.
The graphical plot of the resonance zone in the
$\varepsilon\Omega$-plane form the so-called {\sl Arnold's tongue}.

By a method similar to Arnold's, one may also estimate the width of  ${\cal
I}_r$ for the off-center reflections
$R_{r,\Omega}$, $0\leq
r\leq 1/3$.  We will show the different behaviors of
$R=R_{r,\pi}$ and
$\overline{R}=\overline{R}_{r,\pi}$ at the same time.

For simplicity of computation, let us first consider the resonance zone containing
$\pi$.  Writing $\Omega = \pi + a$, the
second iterates of the maps are
\begin{align*}
R^2(x) &= x + 2\pi + 2a - 2\sum_{k=1}^\infty \frac{r^k}{k}\sin(kx) -
2\sum_{k=1}^\infty \frac{r^k}{k}\sin k\left(x+\pi+a-2\sum
\frac{r^k}{k}\sin(kx)\right) \\
\overline{R}^2(x) &= x + 2\pi + 2a + 2\sum_{k=1}^\infty \frac{r^k}{k}\sin(kx) +
2\sum_{k=1}^\infty \frac{r^k}{k}\sin k\left(x+\pi+a+2\sum
\frac{r^k}{k}\sin(kx)\right)
\end{align*}
The equations of resonance are $R^2(x) = x + 2\pi$ and
$\overline{R}^2(x) = x + 2\pi$.
Let $v = a \mp 2\sum \frac{r^k}{k}\sin(kx)$, we have
$$
0 = v \pm \sum_{k=1}^\infty \left[ \frac{r^k}{k}\sin(kx) -
\frac{(-r)^k}{k}\sin k(x+v)\right],
$$
where
$v = v_1r + v_2r^2 + v_3r^3 + \cdots.$
Note that the solutions of $v$'s for $R$ and $\overline{R}$
do not only differ by a sign.  One can see this by the subtle
combinations of the signs of the infinite series in their second iterates. 
Inductively, one may show that for
$\overline{R}$, we have
$v_k = \dfrac{2}{k}\sin(kx)$, while those for $R$ are
\begin{align*}
v_1 &= -2\sin(x), \\
v_2 &= \sin(2x), \\
v_3 &= 2\sin(x) - \frac{8}{3}\sin(3x), \\
\end{align*}
The above result leads to a $r$-series for $a$ and its maximum and
minimum provide bounds for the resonance zone, namely,
$$
a = \begin{cases}
2\sin(2x)r^2 + \left[2\sin(x) - \frac{7}{3}\sin(3x)\right]r^3 + \cdots,
&\text{for the map $R$,}\\
0 &\text{for the map $\overline{R}$.}
\end{cases}
$$
This calculation agrees with our previous remark that mode-locking
(near $\omega = \pi$) does not occur for
$\overline{R}$.  Furthermore, we have
\begin{theorem}\label{ThmOCwidth}
The width of ${\cal I}_{r}$ is bounded by $Cr^2$
for $\Omega=\pi$ and $Cr$ for general $\Omega$.
\end{theorem}
The computation for the resonance zone at a general $\Omega=2p\pi/q$ is more
complicated.  The equation to formally expand is
$R_{r,a+2p\pi/q}^q(x) = x + 2p\pi$.  The coefficients $a_k$ of
$$
a = a_1r + a_2r^2 + a_3r^3 + \cdots
$$
provide the estimates of ${\cal I}_r$.  It turns out that the first term
$a_1$ does not vanish, indeed,
$$
q a_1 = 2 \sum_{j=0}^{q-1} \sin\left(x+\frac{2jp\pi}{q}\right).
$$
This may not be a sharp estimate, yet we can only conclude that the width
of ${\cal I}_r$ is of order $r$ in general.

\clearpage

\end{document}